\documentclass[]{interact}
\usepackage{lmodern}

\usepackage{epstopdf}
\usepackage{amsfonts,amsbsy,amscd,amsgen,amsthm,dsfont,amsmath,amssymb,mathrsfs,mathtools,array}
\usepackage{url,bm,multirow,esvect,cancel,enumitem}
\usepackage[caption=false]{subfig}

\usepackage[numbers,sort&compress]{natbib}
\bibpunct[, ]{[}{]}{,}{n}{,}{,}
\renewcommand\bibfont{\fontsize{10}{12}\selectfont}
\makeatletter
\def\NAT@def@citea{\def\@citea{\NAT@separator}}
\makeatother

\theoremstyle{plain}
\newtheorem{theorem}{Theorem}[section]
\newtheorem{lemma}[theorem]{Lemma}
\newtheorem{corollary}[theorem]{Corollary}

\theoremstyle{definition}

\newtheorem{example}[theorem]{Example}

\theoremstyle{remark}
\newtheorem{remark}{Remark}

\begin{document}

\articletype{ORIGINAL PAPER}

\title{New formulas for moments of the multivariate normal distribution extending Stein's lemma and Isserlis theorem}

\author{
\name{Konstantinos Mamis\textsuperscript{a}\thanks{Email: kmamis@uw.edu}}
\affil{\textsuperscript{a}Department of Applied Mathematics, University of Washington, Seattle, WA 98195, USA}
}

\maketitle

\begin{abstract}
We prove a formula for the evaluation of expectations containing a scalar function of a Gaussian random vector multiplied by a product of the random vector components, each one raised to a non-negative integer power. Some of the powers could be of zeroth order, and, for expectations containing only one vector component to the first power,  the formula reduces to Stein's lemma for the multivariate normal distribution. Furthermore, by setting the function inside expectation equal to one, we easily re-derive Isserlis theorem and its generalizations,  regarding higher-order moments of a Gaussian random vector. We provide two proofs of the formula, the first being a rigorous proof via mathematical induction.  The second proof is a formal, constructive derivation based on treating the expectation not as an integral, but as the consecutive actions of pseudodifferential operators defined via the moment-generating function of the Gaussian random vector.
\end{abstract}

\begin{keywords}
normal distribution; Stein's lemma; Isserlis theorem; Hermite polynomials; pseudodifferential operator; multi-index; index-matrix\end{keywords}

\section{Introduction}\label{sec:intro}
Stein's lemma \cite{stein}  is a well-known identity of the normal distribution, with applications in statistics, see e.g. \cite{brown}, \cite{zhang}, \cite{new_stein}. More specifically, it constitutes the starting point for the celebrated Stein's method on the distance between distributions, see e.g. \cite{chen2011}. For the case of a scalar random variable $X$ that follows the univariate normal distribution $\mathcal{N}(\mu,\sigma^2)$, Stein's lemma reads
\begin{equation}\label{eq:stein_uni}
\mathsf{E}\left[g(X)X\right]=\mu\mathsf{E}[g(X)]+\sigma^2\mathsf{E}\left[g'(X)\right],
\end{equation}
where $\mathsf{E}\left[\bm{\cdot}\right]$ is the expectation operator, and prime denotes the first derivative of function $g$. In our recent work \cite{mamis_stein}, we extended scalar Stein's lemma \eqref{eq:stein_uni} for expectations containing $X$ at an integer power $n$:
\begin{equation}\label{eq:stein_uni_new}
\mathsf{E}\left[g(X)X^n\right]=\sum_{\ell=0}^n\binom{n}{\ell}\mu^{n-\ell}\sum_{k=0}^{\lfloor\ell/2\rfloor}H_{\ell,k}\sigma^{2(\ell-k)}\mathsf{E}\left[g^{(\ell-2k)}(X)\right], \ \ n\in\mathbb{N},
\end{equation}
where $g^{(\ell)}$ denotes the $\ell$th derivative of $g$, $\binom{n}{\ell}=\frac{n!}{\ell!(n-\ell)!}$ is the binomial coefficient, $\lfloor\bm{\cdot}\rfloor$ is the floor function, and the numbers 
\begin{equation}\label{eq:hermite_numbers}
H_{\ell,k}=\frac{\ell!}{2^kk!(\ell-2k)!}, \ \ k=0,\ldots,\lfloor \ell/2\rfloor,
\end{equation}
are the \textit{signless Hermite coefficients}, i.e., the absolute values of the coefficients appearing in the probabilist's Hermite polynomial of $\ell$th order \cite[expression 22.3.11]{stegun},
$\mathrm{He}_{\ell}(x)=\sum_{k=0}^{\lfloor \ell/2\rfloor}(-1)^k H_{\ell,k}x^{\ell-2k}$. For applications of Hermite polynomials in probability theory, see \cite[Sec. 2.11.5]{shiryaev}. Results similar to Eq.~\eqref{eq:stein_uni_new} have also been derived in \cite{Jamol1,Jamol2}, using the Rodrigues formula. Other generalizations of Stein's lemma \eqref{eq:stein_uni} involving Hermite polynomials have been given by \cite{Goldstein}, and more generally by \cite{Azmoodeh}.

In \cite{mamis_stein}, we proved the scalar generalized Stein's lemma \eqref{eq:stein_uni_new} by mathematical induction on the index $n$. In order to provide insight on the derivation of Eq.~\eqref{eq:stein_uni_new}, we additionally presented a constructive proof; this was based on the alternative definition of the mean value of a deterministic function $g(X)$ with scalar Gaussian argument $X\sim \mathcal{N}(\mu,\sigma^2)$ as the action of a pseudodifferential operator:
\begin{equation}\label{eq:uni_pseudo}
    \mathsf{E}[g(X)]=\left.\exp\left(\frac{\sigma^2}{2}\frac{d^2}{dx^2}\right)g(x)\right|_{x=\mu}.
\end{equation}
In the present work, we extend both proof techniques to the multidimensional case, in order to derive the multivariate counterpart of Eq.~\eqref{eq:stein_uni_new}, that is, a formula for $\mathsf{E}\left[g(\bm{X})\prod_{i=1}^NX_i^{n_i}\right]$, where $\bm{X}$ is an $N$-dimensional Gaussian random vector $\bm{X}\sim \mathcal{N}(\bm{\mu},\bm{C})$. To the best of our knowledge, such a formula has not been derived before, and  constitutes the generalization of the multivariate Stein's lemma \cite{liu}:
\begin{equation}\label{eq:stein}
\mathsf{E}\left[g(\bm{X})X_i\right]=\mu_i\mathsf{E}\left[g(\bm{X})\right]+\sum_{j=1}^NC_{i j}\mathsf{E}\left[\partial_j g(\bm{X})\right], \ \ i=1,\ldots,N,
\end{equation}
with $\partial_jg(\bm{X})=\partial g(\bm{X})/\partial X_j$.  In our generalization of Eq.~\eqref{eq:stein}, the Gaussianity of the random vector $\bm{X}$ is retained; we have to note that a formula has recently been derived that evaluates $\mathsf{E}[g(\bm{X})X_1^2]$ for $\bm{X}$ following a multivariate elliptical distribution \cite{zuo}. 

The derivation and proof of the generalization of Stein's lemma \eqref{eq:stein} relies heavily on the use of multi-index and index-matrix notation, introduced in Sec.~\ref{sec:multiindex}. Our main result, the multivariate generalized Stein's lemma, is presented in Sec.~\ref{sec:result}. In Sec.~\ref{sec:corollaries}, we re-derive known results for moments of Gaussian vectors, such as Isserlis theorem, as corollaries to the generalized Stein's lemma. Rigorous proof of generalized Stein's lemma via multidimensional mathematical induction is presented in Sec.~\ref{sec:induction}. In Sec.~\ref{sec:construction}, we generalize the pseudodifferential operator relation \eqref{eq:uni_pseudo} for expectations involving Gaussian vectors, and present an alternative, constructive derivation for our main result. This use of pseudodifferential operators in the multidimensional case constitutes the main methodological novelty of the present work.  In the concluding discussion of Sec.~\ref{sec:conclusions}, we discuss the use of infinite-dimensional counterparts of Stein's lemma and its generalization in deriving evolution equations for the probability density function of the response to dynamical systems under Gaussian random excitation. 

\subsection{Multi-index and index-matrix notation}\label{sec:multiindex}
Expressing and proving the generalization of  Stein's lemma in $N$ dimensions is simplified by using the multi-index \cite[p. 319]{reed} and index-matrix \cite{Rahman2017, Slepian1972} notation. A multi-index $\bm{n}$ and an index-matrix $\bm{L}$ are an $N$-dimensional index vector and an $N\times N$ matrix of indices respectively
\[\bm{n}=\left\{n_i\right\}_{i=1}^N\in\mathbb{N}^N_0, \ \ \bm{L}=\left\{\ell_{ij}\right\}_{i,j=1}^N\in\mathbb{N}^{N\times N}_0.\]
For multi-indices and index-matrices, we define the following operations:
\begin{itemize}[leftmargin=*]
    \item Partial ordering of multi-indices
    \[\bm{m}\leq\bm{n}\Leftrightarrow m_i\leq n_i, \ \ \forall i\in\{1,\ldots,N\}.\]
    \item Partial ordering of index-matrices
    \[\bm{K}\leq\bm{L}\Leftrightarrow k_{ij}\leq\ell_{ij}, \ \ \forall i,j\in\{1,\ldots,N\}.\]
    \item Factorial of a multi-index \[\bm{n}!=\prod_{i=1}^Nn_i!.\]
    \item Factorial of an index-matrix \[\bm{L}!=\prod_{i=1}^N\prod_{j=1}^N\ell_{ij}!.\]
    \item  Row-sum $\bm{r}(\bm{L})$ and column-sum $\bm{c}(\bm{L})$ vectors of index-matrix $\bm{L}$
\begin{equation}\label{eq:row_column_sum}
r_i(\bm{L})=\sum_{j=1}^N\ell_{ij},\ \  c_i(\bm{L})=\sum_{j=1}^N\ell_{ji}, \ \  i=1,\ldots,N.
\end{equation}
    \item Multinomial coefficient of a multi-index and an index-matrix with respect to the rows of the index-matrix ($r$-multinomial coefficient)
    \begin{equation}\label{eq:multi_r}
\binom{\bm{n}}{\bm{L}}_r=\frac{\bm{n}!}{[\bm{n}-\bm{r}(\bm{L})]!\bm{L}!}.
\end{equation}
    \item Raising an $N$-dimensional vector $\bm{X}=\{X_i\}_{i=1}^N$ to a multi-index power $\bm{n}$ \[\bm{X}^{\bm{n}}=\prod_{i=1}^NX_i^{n_i}.\]
    \item Raising an $N\times N$ matrix $\bm{C}=\{C_{ij}\}_{i,j=1}^N$ to an index-matrix power $\bm{L}$
    \[\bm{C}^{\bm{L}}=\prod_{i=1}^N\prod_{j=1}^NC_{ij}^{\ell_{ij}}.\]
    \item Partial derivative of multi-index order $\bm{n}$ of a scalar function $g(\bm{X})$ with $N$-dimensional vector argument $\bm{X}=\{X_i\}_{i=1}^N$
    \[\partial^{\bm{n}}g(\bm{X})=\prod_{i=1}^N\partial_i^{n_i}g(\bm{X}).\]
\end{itemize}
\section{Results}
Our main result is the generalization of multivariate Stein's lemma \eqref{eq:stein}, stated in the following Theorem.
\subsection{Multivariate generalized Stein’s lemma}\label{sec:result}
\begin{theorem}\label{cor:index_matrix}
For $\bm{X}\sim\mathcal{N}(\bm{\mu},\bm{C})$, a smooth enough function $g:\mathbb{R}^N\rightarrow\mathbb{R}$, a multi-index $\bm{n}\in\mathbb{N}_0^N$, and under the assumption that all expectations involved exist, it holds true that
\begin{align}\label{eq:stein_new_new}
\mathsf{E}\left[g(\bm{X})\bm{X}^{\bm{n}}\right]=&\sum_{\substack{\bm{L}\in\mathbb{N}^{N\times N}_0\\\bm{r}(\bm{L})\leq\bm{n}}} \binom{\bm{n}}{\bm{L}}_r\bm{\mu}^{\bm{n}-\bm{r}(\bm{L})}\times\nonumber\\&\times\sum_{\substack{\bm{K}\in\mathbb{N}_{0,\mathrm{sym}}^{N\times N}\\\bm{K}\leq\bm{K}^{\max}(\bm{L})}}H_{\bm{L},\bm{K}}\bm{C}_U^{\bm{L}+\bm{L}^T-\mathrm{diag}(\bm{L})-\bm{K}}\mathsf{E}\left[\partial^{\bm{c}(\bm{L}-\bm{K}-\mathrm{diag}(\bm{K}))}g(\bm{X})\right].
\end{align}
Subscript $U$ denotes the upper triangular part of a matrix, superscript $T$ denotes the matrix transpose, and $\mathrm{diag}(\cdot)$ denotes the diagonal part of a matrix. The inner sum in the right-hand side of Eq.~\eqref{eq:stein_new_new} is over all \textbf{symmetric} index-matrices $\bm{K}$, with  $\bm{K}^{\max}(\bm{L})$ defined as
\begin{equation}\label{eq:K_max}
k^{\max}_{ii}=\lfloor\ell_{ii}/2\rfloor, \ \ \text{and} \ \ k^{\max}_{ij}=\min\{\ell_{ij},\ell_{ji}\} \ \ \text{for} \ \ i\neq j.
\end{equation}
Coefficients $H_{\bm{L},\bm{K}}$, that constitute the multidimensional counterpart of signless Hermite coefficients \eqref{eq:hermite_numbers}, are defined as
\begin{equation}\label{eq:DD}
H_{\bm{L},\bm{K}}=\frac{\bm{L}!}{2^{\mathrm{tr}(\bm{K})}\bm{K}_U!(\bm{L}-{\bm{K}}-\mathrm{diag}(\bm{K}))!},
\end{equation} 
where $\mathrm{tr}(\bm{K})=\sum_{i=1}^Nk_{ii}$ is the trace of index-matrix $\bm{K}$. 
\end{theorem}
In order to have less lengthy expressions in the process of proving Eq.~\eqref{eq:stein_new_new}, we define the symmetric index-matrices
$\tilde{\bm{L}}=\left\{\tilde{\ell}_{ij}\right\}_{i,j=1}^N$, $\tilde{\bm{K}}=\left\{\tilde{k}_{ij}\right\}_{i,j=1}^N$:
\begin{equation}\label{eq:L_tilde}
\tilde{\bm{L}}=\bm{L}+\bm{L}^T-\mathrm{diag}(\bm{L})\Rightarrow\tilde{\ell}_{ii}=\ell_{ii}, \ \ \text{and} \ \ \tilde{\ell}_{ij}=\ell_{ij}+\ell_{ji} \ \ \text{for} \ \ i\neq j,
\end{equation}
\begin{equation}\label{eq:K_tilde}
\tilde{\bm{K}}=\bm{K}+\mathrm{diag}(\bm{K})\Rightarrow\tilde{k}_{ii}=2k_{ii}, \ \ \text{and} \ \ \tilde{k}_{ij}=k_{ij} \ \ \text{for} \ \ i\neq j.   
\end{equation}
Using $\tilde{\bm{L}}$, $\tilde{\bm{K}}$, multivariate generalized Stein's lemma \eqref{eq:stein_new_new} is equivalently expressed as
\begin{align}\label{eq:stein_new}
\mathsf{E}\left[g(\bm{X})\bm{X}^{\bm{n}}\right]=&\sum_{\substack{\bm{L}\in\mathbb{N}^{N\times N}_0\\\bm{r}(\bm{L})\leq\bm{n}}} \binom{\bm{n}}{\bm{L}}_r\bm{\mu}^{\bm{n}-\bm{r}(\bm{L})}\sum_{\substack{\bm{K}\in\mathbb{N}_{0,\mathrm{sym}}^{N\times N}\\\bm{K}\leq\bm{K}^{\max}(\bm{L})}}H_{\bm{L},\bm{K}}\bm{C}_U^{\tilde{\bm{L}}-\bm{K}}\mathsf{E}\left[\partial^{\bm{c}(\bm{L}-\tilde{\bm{K}})}g(\bm{X})\right],
\end{align}
with
\begin{equation}\label{eq:D}
H_{\bm{L},\bm{K}}=\frac{\bm{L}!}{2^{\mathrm{tr}(\bm{K})}\bm{K}_U!(\bm{L}-\tilde{\bm{K}})!}.
\end{equation} 
We illustrate the application of Eq.~\eqref{eq:stein_new} in the following example.
\begin{example}\label{example}
For $N=2$, $\bm{n}=(1,2)$, $\bm{\mu}=\bm{0}$, and with $\sigma^2_1:=C_{11}$, $\sigma^2_2:=C_{22}$,  Eq.~\eqref{eq:stein_new} results in
\begin{align}\label{eq:example}
\mathsf{E}\left[g(X_1,X_2)X_1X_2^2\right]=\left(\sigma_1^2\sigma_2^2+2C_{12}^2\right)\mathsf{E}\left[\frac{\partial g(X_1,X_2)}{\partial X_1}\right]+3\sigma_2^2C_{12}\mathsf{E}\left[\frac{\partial g(X_1,X_2)}{\partial X_2}\right]&+\nonumber\\+\left(2\sigma_1^2\sigma_2^2C_{12}+C_{12}^3\right)\mathsf{E}\left[\frac{\partial^3 g(X_1,X_2)}{\partial X_1^2\partial X_2}\right]+\left(\sigma_1^2\sigma_2^4+2\sigma_2^2C_{12}^2\right)\mathsf{E}\left[\frac{\partial^3 g(X_1,X_2)}{\partial X_1\partial X_2^2}\right]&+\nonumber\\+\sigma_1^2C_{12}^2\mathsf{E}\left[\frac{\partial^3 g(X_1,X_2)}{\partial X_1^3}\right]+\sigma_2^4C_{12}\mathsf{E}\left[\frac{\partial^3 g(X_1,X_2)}{\partial X_2^3}\right]&.
\end{align} 
Since $\bm{\mu}=\bm{0}$, the outer sum in the right-hand side of Eq.~\eqref{eq:stein_new} extends over all index-matrices $\bm{L}\in\mathbb{N}_0^{N\times N}$ with $\bm{r}(\bm{L})=\bm{n}=(1,2)$. These matrices are $\begin{pmatrix}1&0\\2&0\end{pmatrix}$, $\begin{pmatrix}1&0\\0&2\end{pmatrix}$, $\begin{pmatrix}1&0\\1&1\end{pmatrix}$, $\begin{pmatrix}0&1\\2&0\end{pmatrix}$, $\begin{pmatrix}0&1\\0&2\end{pmatrix}$ and $\begin{pmatrix}0&1\\1&1\end{pmatrix}$.  For each of these index-matrices,  we calculate, in Table \ref{table:1}, the rest of the quantities appearing in each term of the $\bm{L}$-sum in the right-hand side of Eq.~\eqref{eq:stein_new}.  By substituting the quantities from Table \ref{table:1} into Eq.~\eqref{eq:stein_new}, and after some regrouping of terms, we obtain Eq.~\eqref{eq:example}. Note that Eq.~\eqref{eq:example} can also be validated by repetitive applications of Stein's lemma, Eq.~\eqref{eq:stein}, for $\mathsf{E}\left[g(X_1,X_2)X_1X_2^2\right]$.
\end{example}
\begin{table}
\footnotesize
\caption{Matrix calculations for Example \ref{example}.}\label{table:1}
\begin{center}
\begin{tabular}{|c|c|c|c|c|c|c|c|} 
 \hline
$\bm{L}$ & $\binom{\bm{n}}{\bm{L}}_r$ & $\tilde{\bm{L}}$& $\bm{K}^{\max}$ & $\bm{K}$ & $\tilde{\bm{K}}$ & $H_{\bm{L},\bm{K}}$ & $\bm{c}(\bm{L}-\tilde{\bm{K}})$\\
\hline\hline
$\begin{pmatrix}1 & 0\\2&0\end{pmatrix}$ & 1 &$\begin{pmatrix}1&2\\2&0\end{pmatrix}$& $\begin{pmatrix}0 & 0\\0 &0\end{pmatrix}$ &$\begin{pmatrix}0 & 0\\0 &0\end{pmatrix}$ &$\begin{pmatrix}0 & 0\\0 &0\end{pmatrix}$ &1&$\begin{pmatrix} 3\\0\end{pmatrix}$\\
\hline
\multirow{2}{*}{$\begin{pmatrix}1&0\\0&2\end{pmatrix}$} &\multirow{2}{*}{1}&\multirow{2}{*}{$\begin{pmatrix}1&0\\0&2\end{pmatrix}$} &\multirow{2}{*}{$\begin{pmatrix}0&0\\0&1\end{pmatrix}$} & $\begin{pmatrix}0&0\\0&0\end{pmatrix}$ & $\begin{pmatrix}0&0\\0&0\end{pmatrix}$ & 1 & $\begin{pmatrix}1\\2\end{pmatrix}$\\
\cline{5-8}
&&&&$\begin{pmatrix}0&0\\0&1\end{pmatrix}$&$\begin{pmatrix}0&0\\0&2\end{pmatrix}$&1&$\begin{pmatrix}1\\0\end{pmatrix}$\\
\hline
$\begin{pmatrix}1&0\\1&1\end{pmatrix}$ & 2&$\begin{pmatrix}1&1\\1&1\end{pmatrix}$&$\begin{pmatrix}0&0\\0&0\end{pmatrix}$&$\begin{pmatrix}0&0\\0&0\end{pmatrix}$&$\begin{pmatrix}0&0\\0&0\end{pmatrix}$&1&$\begin{pmatrix}2\\1\end{pmatrix}$\\
\hline
\multirow{2}{*}{$\begin{pmatrix}0&1\\2&0\end{pmatrix}$}&\multirow{2}{*}{1}&\multirow{2}{*}{$\begin{pmatrix}0&3\\3&0\end{pmatrix}$}&\multirow{2}{*}{$\begin{pmatrix}0&1\\1&0\end{pmatrix}$}&$\begin{pmatrix}0&0\\0&0\end{pmatrix}$&$\begin{pmatrix}0&0\\0&0\end{pmatrix}$&1&$\begin{pmatrix}2\\1\end{pmatrix}$\\
\cline{5-8}
&&&&$\begin{pmatrix}0&1\\1&0\end{pmatrix}$&$\begin{pmatrix}0&1\\1&0\end{pmatrix}$&2&$\begin{pmatrix}1\\0\end{pmatrix}$\\
\hline
\multirow{2}{*}{$\begin{pmatrix}0&1\\0&2\end{pmatrix}$}&\multirow{2}{*}{1}&\multirow{2}{*}{$\begin{pmatrix}0&1\\1&2\end{pmatrix}$}&\multirow{2}{*}{$\begin{pmatrix}0&0\\0&1\end{pmatrix}$}&$\begin{pmatrix}0&0\\0&0\end{pmatrix}$&$\begin{pmatrix}0&0\\0&0\end{pmatrix}$&1&$\begin{pmatrix}0\\3\end{pmatrix}$\\
\cline{5-8}
&&&&$\begin{pmatrix}0&0\\0&1\end{pmatrix}$&$\begin{pmatrix}0&0\\0&2\end{pmatrix}$&1&$\begin{pmatrix}0\\1\end{pmatrix}$\\
\hline
\multirow{2}{*}{$\begin{pmatrix}0&1\\1&1\end{pmatrix}$}&\multirow{2}{*}{2}&\multirow{2}{*}{$\begin{pmatrix}0&2\\2&1\end{pmatrix}$}&\multirow{2}{*}{$\begin{pmatrix}0&1\\1&0\end{pmatrix}$} &$\begin{pmatrix}0&0\\0&0\end{pmatrix}$&$\begin{pmatrix}0&0\\0&0\end{pmatrix}$ &1&$\begin{pmatrix}1\\2\end{pmatrix}$\\
\cline{5-8}
&&&&$\begin{pmatrix}0&1\\1&0\end{pmatrix}$&$\begin{pmatrix}0&1\\1&0\end{pmatrix}$&1&$\begin{pmatrix}0\\1\end{pmatrix}$\\
\hline
\end{tabular}
\end{center}
\end{table}
In Sec.~\ref{sec:induction}, we prove Eq.~\eqref{eq:stein_new} rigorously, by multidimensional mathematical induction on $\bm{n}\in\mathbb{N}^N_0$. In addition to this proof, and in order to provide the reader with more insight on the derivation of Eq.~\eqref{eq:stein_new},  we present, in Sec.~\ref{sec:construction}, a constructive formal proof of Theorem \ref{cor:index_matrix}. This constructive proof is based on treating the mean value operator not as an integral, but as the sequential action of a number of pseudodifferential operators that are introduced in Theorem \ref{def1} via the moment-generating function of the Gaussian random vector.  The action of these pseudodifferential operators is determined by their Taylor series expansions, under the formal assumption that all infinite series involved are summable. 
\subsection{Corollaries for moments of a Gaussian vector}\label{sec:corollaries}
From Theorem \ref{cor:index_matrix} we re-derive existing results regarding moments of Gaussian vectors, such as Isserlis theorem.
\begin{corollary}\label{cor:song}
For the Gaussian random vector $\bm{X}\sim\mathcal{N}(\bm{\mu},\bm{C})$, the following formula for its product moments holds true:
\begin{equation}\label{eq:song}
\mathsf{E}\left[\bm{X}^{\bm{n}}\right]=\sum_{\substack{\bm{K}\in\mathbb{N}_{0,\mathrm{sym}}^{N\times N}\\\bm{r}(\tilde{\bm{K}})\leq\bm{n}}}d_{\bm{n},\bm{K}}\bm{\mu}^{\bm{n}-\bm{r}(\tilde{\bm{K}})}\bm{C}_U^{\bm{K}},
\end{equation}
where $\tilde{\bm{K}}$ is defined from $\bm{K}$ using Eq.~\eqref{eq:K_tilde}, and
\begin{equation}\label{d_coeff}
d_{\bm{n},\bm{K}}=\frac{\bm{n}!}{2^{\mathrm{tr}(\bm{K})}\bm{K}_U![\bm{n}-\bm{r}(\tilde{\bm{K}})]!}.
\end{equation}
Summation in the right-hand side of Eq.~\eqref{eq:song} extends over all index-matrices $\bm{K}$ for which the respective index-matrix $\tilde{\bm{K}}$ satisfies the condition $\bm{r}(\tilde{\bm{K}})\leq\bm{n}$. 
\end{corollary}
\begin{proof}
We derive Eq.~\eqref{eq:song} by setting $g(\bm{X})=1$ in formula \eqref{eq:stein_new}. By setting $g(\bm{X})=1$, all derivatives in the right-hand side of Eq.~\eqref{eq:stein_new} are zero, except for the zeroth order derivative, that is for $\bm{c}(\bm{L}-\tilde{\bm{K}})=\bm{0}$. By the definition relation \eqref{eq:K_tilde} of $\tilde{\bm{K}}$, $\bm{c}(\bm{L}-\tilde{\bm{K}})=\bm{0}$ is achieved for $\bm{L}=\tilde{\bm{K}}$. Last, coefficients $d_{\bm{n},\bm{K}}$ are calculated as
\[\binom{\bm{n}}{\tilde{\bm{K}}}_rH_{\tilde{\bm{K}},\bm{K}}=\frac{\bm{n}!}{[\bm{n}-\bm{r}(\tilde{\bm{K}})]!\tilde{\bm{K}}!}\frac{\tilde{\bm{K}}!}{2^{\mathrm{tr}(\bm{K})}\bm{K}_U!}=d_{\bm{n},\bm{K}}.\]
Substitution of the above in Eq.~\eqref{eq:stein_new} results in Eq.~\eqref{eq:song}.
\end{proof}
Eq.~\eqref{eq:song} has been proven by Song and Lee in \cite{song}, using Price's theorem, see \cite{price}, \cite[Sec. 5.3.2]{song2}; see also the recent work \cite{ken}.

\begin{corollary}\label{cor:isserlis}
From Eq.~\eqref{eq:song}, we derive the formula for the higher order moments of an $N$-dimensional Gaussian random vector $\bm{X}$ with zero mean value:
\begin{equation}\label{eq:isserlis}
\mathsf{E}\left[\prod_{i=1}^NX_i\right]= \left\{
\begin{array}{ll}
      0 & \text{for } N\text{ odd}, \\
      \sum_{P\in\wp_N^2}\prod_{(i,j)\in P}C_{ij} & \text{for } N\text{ even},\\
\end{array} 
\right. 
\end{equation}
with $\wp_N^2$ being the set of all partitions of $\{1,\ldots,N\}$ into unordered pairs.
\end{corollary}
\begin{proof}
For $\bm{\mu}=\bm{0}$ and $\bm{n}=\bm{1}$, the summation in Eq.~\eqref{eq:song} extends over all $\bm{K}$ with $\bm{r}(\tilde{\bm{K}})=\bm{1}$. Since the diagonal elements of $\tilde{\bm{K}}$ are even numbers (see definition relation \eqref{eq:K_tilde}), $\tilde{\bm{K}}$ in this case has zero diagonal elements, and it is equal to $\bm{K}$. In addition, we calculate $d_{\bm{1},\bm{K}}=1$. Thus, Eq.~\eqref{eq:song} reads
\begin{equation}\label{eq:is}
\mathsf{E}\left[\prod_{i=1}^NX_i\right]=\sum_{\bm{K}\in\mathscr{A}_N}\bm{C}_U^{\bm{K}},
\end{equation}
where $\mathscr{A}_N$ is the set of all $N\times N$ matrices that are i) symmetric, ii) have all diagonal elements zero, iii) their elements are either 0 or 1, iv) each row sum equals to one. Thus, matrices $\bm{K}\in\mathscr{A}_N$ are identified \cite[definition 2.1]{biggs} as the adjacency matrices of undirected 1-regular graphs between $N$ nodes. By virtue of the handshaking lemma \cite[Theorem 2.1]{graph}, the number of nodes $N$ cannot be of the same parity with 1, which is the graph degree. Thus, for $N$ odd, the set $\mathscr{A}_N$ is empty and so the odd moments are zero. For even $N$, each undirected 1-regular graph between $N$ nodes is equivalent to one partition of set $\{1,\ldots,N\}$ into unordered pairs, and so Eq.~\eqref{eq:is} is expressed equivalently as the branch of Eq.~\eqref{eq:isserlis} for even $N$.
\end{proof}
Eq.~\eqref{eq:isserlis} is the Isserlis theorem \cite{isserlis},  also known in physics literature as Wick's theorem \cite{wick}.  For a review of the related literature, see also \cite{triant}. 

\section{Proof of Theorem \ref{cor:index_matrix} by mathematical induction}\label{sec:induction}
We prove Theorem \ref{cor:index_matrix} by multidimensional mathematical induction on the multi-index of  exponents $\bm{n}$. For this proof, we introduce the multi-index $\bm{e}^{(i)}$ which has all its components equal to zero, except for the $i$th component which is equal to one. Similarly, we introduce the index-matrix $\bm{E}^{(ij)}$ that has all its elements equal to zero, except for the $ij$ element which is equal to one. Last, index-matrix $\bm{E}^{(ij)}_{\text{sym}}$ is symmetric, having both $ij$ and $ji$ elements equal to one, and the rest of its elements are equal to zero. 
\begin{itemize}[leftmargin=*]
    \item Base case: $\bm{n}=\bm{e}^{(i)}$, $i=1,\ldots,N$.
\end{itemize}
The index-matrices $\bm{L}$ with $\bm{r}(\bm{L})\leq\bm{e}^{(i)}$ are: i) the zero matrix $\bm{0}$, ii) the matrices $\bm{E}^{(ij)}$, $j=1,\ldots,N$. Matrix $\bm{0}$ results in the term $\mu_i\mathsf{E}\left[g(\bm{X})\right]$ in the $\bm{L}$-sum of Eq.~\eqref{eq:stein_new}. For each matrix $\bm{E}^{(ij)}$, $j=1,\ldots,N$, we calculate that $\binom{\bm{e}^{(i)}}{\bm{E}^{(ij)}}_r=1$, $\bm{K}^{\max}(\bm{E}^{(ij)})=\bm{K}=\tilde{\bm{K}}=\bm{0}$, $H_{\bm{E}^{(ij)},\bm{0}}=1$, $\tilde{\bm{L}}=\bm{E}^{(ij)}_{\text{sym}}$, $\bm{C}_U^{\bm{E}^{(ij)}_{\text{sym}}}=C_{ij}$, and $\bm{c}(\bm{E}^{(ij)})=\bm{e}^{(j)}$. Thus, each index-matrix $\bm{E}^{(ij)}$ results in the term $C_{ij}\mathsf{E}\left[\partial_jg(\bm{X})\right]$ in the $\bm{L}$-sum. Summation of all the said terms results in multivariate Stein's lemma, Eq.~\eqref{eq:stein}. Eq.~\eqref{eq:stein} is easily proven
by integration by parts, using the properties of the derivatives of the $N$-variate normal distribution, see \cite[Eq.347]{cookbook}.
\begin{itemize}[leftmargin=*]
\item Inductive hypothesis: Eq.~\eqref{eq:stein_new} holds true for $\bm{n}$.
\end{itemize}
\begin{itemize}[leftmargin=*]
\item Inductive step: Prove that Eq.~\eqref{eq:stein_new} holds true for $\bm{n}+\bm{e}^{(i)}$, $i=1,\ldots,N$.
\end{itemize}
By using the inductive hypothesis, we have
\begin{align}\label{eq:proof1}
&\mathsf{E}\left[g(\bm{X})\bm{X}^{\bm{n}+\bm{e}^{(i)}}\right]=\mathsf{E}\left[(g(\bm{X})X_i)\bm{X}^{\bm{n}}\right]=\nonumber\\&=\sum_{\bm{r}(\bm{L})\leq\bm{n}} \binom{\bm{n}}{\bm{L}}_r\bm{\mu}^{\bm{n}-\bm{r}(\bm{L})}\times\sum_{\bm{K}\leq\bm{K}^{\max}(\bm{L})}H_{\bm{L},\bm{K}}\bm{C}_U^{\tilde{\bm{L}}-\bm{K}}\mathsf{E}\left[\partial^{\bm{c}(\bm{L}-\tilde{\bm{K}})}(g(\bm{X})X_i)\right].
\end{align}
By the general Leibniz rule \cite[expression 3.3.8]{stegun}, we calculate the derivative 
\begin{equation}\label{eq:leibniz}
\partial_i^{c_i(\bm{L}-\tilde{\bm{K}})}\left(g(\bm{X})X_i\right)=\sum_{p=0}^{c_i(\bm{L}-\tilde{\bm{K}})}\binom{c_i(\bm{L}-\tilde{\bm{K}})}{p}\left(\partial_i^{c_i(\bm{L}-\tilde{\bm{K}})-p}g(\bm{X})\right)\left(\partial_i^pX_i\right).
\end{equation}
Since $\partial_i^0X_i=X_i$, $\partial_iX_i=1$, and $\partial_i^pX_i=0$ for $p\geq2$, Eq.~\eqref{eq:leibniz} is simplified into
\begin{align}\label{eq:leibniz2}
\partial_i^{c_i(\bm{L}-\tilde{\bm{K}})}\left(g(\bm{X})X_i\right)&=X_i\partial_i^{c_i(\bm{L}-\tilde{\bm{K}})}g(\bm{X})+{c_i(\bm{L}-\tilde{\bm{K}})}\partial_i^{{c_i(\bm{L}-\tilde{\bm{K}})}-1}g(\bm{X})\nonumber\\&=X_i\partial_i^{c_i(\bm{L}-\tilde{\bm{K}})}g(\bm{X})+\sum_{j=1}^N\left[\ell_{ji}-(1+\delta_{ij})k_{ij}\right]\partial_i^{{c_i(\bm{L}-\tilde{\bm{K}})}-1}g(\bm{X}),
\end{align}
where $\delta_{ij}$ is Kronecker's delta. By using Eq.~\eqref{eq:leibniz2}, we rewrite Eq.~\eqref{eq:proof1} as
\begin{equation}\label{eq:AB}
\mathsf{E}\left[(g(\bm{X})X_i)\bm{X}^{\bm{n}}\right]=A+\sum_{j=1}^NB_j,
\end{equation}
with
\begin{equation}\label{eq:A}
A=\sum_{\bm{r}(\bm{L})\leq\bm{n}} \binom{\bm{n}}{\bm{L}}_r\bm{\mu}^{\bm{n}-\bm{r}(\bm{L})}\sum_{\bm{K}\leq\bm{K}^{\max}(\bm{L})}H_{\bm{L},\bm{K}}\bm{C}_U^{\tilde{\bm{L}}-\bm{K}}\mathsf{E}\left[X_i\partial^{\bm{c}(\bm{L}-\tilde{\bm{K}})}g(\bm{X})\right]
\end{equation}
and
\begin{align}\label{eq:B_j}
B_j=\sum_{\bm{r}(\bm{L})\leq\bm{n}} \binom{\bm{n}}{\bm{L}}_r\bm{\mu}^{\bm{n}-\bm{r}(\bm{L})}\sum_{\bm{K}\leq\bm{K}^{\max}(\bm{L})}&\left[\ell_{ji}-(1+\delta_{ij})k_{ij}\right]H_{\bm{L},\bm{K}}\bm{C}_U^{\tilde{\bm{L}}-\bm{K}}\times\nonumber\\&\times\mathsf{E}\left[\partial_i^{c_i(\bm{L}-\tilde{\bm{K}})-1}\prod_{\substack{p=1\\ p\neq i}}^N\partial_p^{c_p(\bm{L}-\tilde{\bm{K}})}g(\bm{X})\right].
\end{align}
In Eq.~\eqref{eq:B_j}, the term in $\bm{K}$-sum is zero for $\ell_{ii}=2k_{ii}$, or $\ell_{ji}=k_{ij}$ for $i\neq j$. In order to exclude zero terms, we update Eq.~\eqref{eq:B_j} to
\begin{align}\label{eq:B_jj}
B_j=\sum_{\bm{r}(\bm{L})\leq\bm{n}} \binom{\bm{n}}{\bm{L}}_r\bm{\mu}^{\bm{n}-\bm{r}(\bm{L})}\sum_{\bm{K}\leq\bm{K}^{\max}(\bm{L}-\bm{E}^{(ji)})}&\left[\ell_{ji}-(1+\delta_{ij})k_{ij}\right]H_{\bm{L},\bm{K}}\bm{C}_U^{\tilde{\bm{L}}-\bm{K}}\times\nonumber\\\times&\mathsf{E}\left[\partial_i^{c_i(\bm{L}-\tilde{\bm{K}})-1}\prod_{\substack{p=1\\ p\neq i}}^N\partial_p^{c_p(\bm{L}-\tilde{\bm{K}})}g(\bm{X})\right].
\end{align}
By performing the change of index-matrix $\bm{K}'=\bm{K}+\bm{E}^{(ij)}_{\text{sym}}$, we recast Eq.~\eqref{eq:B_jj} into
\begin{align}\label{eq:B_j2}
B_j=&\sum_{\bm{r}(\bm{L})\leq\bm{n}} \binom{\bm{n}}{\bm{L}}_r\bm{\mu}^{\bm{n}-\bm{r}(\bm{L})}\sum_{\bm{E}^{(ij)}_{\text{sym}}\leq\bm{K}'\leq\bm{K}^{\max}(\bm{L}-\bm{E}^{(ji)})+\bm{E}^{(ij)}_{\text{sym}}}\left[\ell_{ji}-(1+\delta_{ij})(k_{ij}'-1)\right]\times\nonumber\\\times &H_{\bm{L},\bm{K}'-\bm{E}^{(ij)}_{\text{sym}}}\bm{C}_U^{\tilde{\bm{L}}-\bm{K}'+\bm{E}^{(ij)}_{\text{sym}}}\times\mathsf{E}\left[\partial_j^{c_j(\bm{L}-\tilde{\bm{K}}')+1}\prod_{\substack{p=1\\ p\neq j}}^N\partial_p^{c_{p}(\bm{L}-\tilde{\bm{K}}')}g(\bm{X})\right].
\end{align}
Since $\lfloor(\ell_{ii}-1)/2\rfloor+1=\lfloor(\ell_{ii}+1)/2\rfloor$, and $\min\{\ell_{ij},\ell_{ji}-1\}=\min\{\ell_{ij}+1,\ell_{ji}\}$ for $i\neq j$, it holds true that $\bm{K}^{\max}(\bm{L}-\bm{E}^{(ji)})+\bm{E}^{(ij)}_{\text{sym}}=\bm{K}^{\max}(\bm{L}+\bm{E}^{(ij)})$. Thus, Eq.~\eqref{eq:B_j2} is expressed equivalently as
\begin{align}\label{eq:B_j3}
B_j=&\sum_{\bm{r}(\bm{L})\leq\bm{n}} \binom{\bm{n}}{\bm{L}}_r\bm{\mu}^{\bm{n}-\bm{r}(\bm{L})}\sum_{\bm{E}^{(ij)}_{\text{sym}}\leq\bm{K}\leq\bm{K}^{\max}(\bm{L}+\bm{E}^{(ij)})}\left[\ell_{ji}-(1+\delta_{ij})(k_{ij}-1)\right]\times\nonumber\\&\times H_{\bm{L},\bm{K}-\bm{E}^{(ij)}_{\text{sym}}}\bm{C}_U^{\tilde{\bm{L}}-\bm{K}+\bm{E}^{(ij)}_{\text{sym}}}\mathsf{E}\left[\partial_j^{c_j(\bm{L}-\tilde{\bm{K}})+1}\prod_{\substack{p=1\\ p\neq j}}^N\partial_p^{c_{p}(\bm{L}-\tilde{\bm{K}})}g(\bm{X})\right].
\end{align}
By applying Stein's lemma \eqref{eq:stein} at the expectation appearing in the right-hand side of Eq.~\eqref{eq:A}, we obtain
\begin{equation}\label{eq:A2}
A=A_0+\sum_{j=1}^NA_j,
\end{equation}
with
\begin{equation}\label{eq:A_0}
A_0=\sum_{\bm{r}(\bm{L})\leq\bm{n}} \binom{\bm{n}}{\bm{L}}_r\bm{\mu}^{\bm{n}+\bm{e}^{(i)}-\bm{r}(\bm{L})}\sum_{\bm{K}\leq\bm{K}^{\max}(\bm{L})}H_{\bm{L},\bm{K}}\bm{C}_U^{\tilde{\bm{L}}-\bm{K}}\mathsf{E}\left[\partial^{\bm{c}(\bm{L}-\tilde{\bm{K}})}g(\bm{X})\right]
\end{equation}
and
\begin{align}\label{A_j}
A_j=\sum_{\bm{r}(\bm{L})\leq\bm{n}} \binom{\bm{n}}{\bm{L}}_r\bm{\mu}^{\bm{n}-\bm{r}(\bm{L})}\sum_{\bm{K}\leq\bm{K}^{\max}(\bm{L})}&H_{\bm{L},\bm{K}}\bm{C}_U^{\tilde{\bm{L}}-\bm{K}+\bm{E}^{(ij)}_{\text{sym}}}\times\nonumber\\&\times\mathsf{E}\left[\partial_j^{c_j(\bm{L}-\tilde{\bm{K}})+1}\prod_{\substack{p=1\\ p\neq j}}^N\partial_p^{c_{p}(\bm{L}-\tilde{\bm{K}})}g(\bm{X})\right].
\end{align}
Thus, under Eqs.~\eqref{eq:AB} and \eqref{eq:A2}, the expectation is expressed as
\begin{equation}\label{eq:AB2}
\mathsf{E}\left[(g(\bm{X})X_i)\bm{X}^{\bm{n}}\right]=A_0+\sum_{j=1}^N(A_j+B_j).
\end{equation}
In order to evaluate each $A_j+B_j$ further, we prove the following lemma.

\begin{lemma}\label{lem:H} For $\bm{L}\in\mathbb{N}^{N\times N}_0$, $\bm{K}\in\mathbb{N}^{N\times N}_{0,\mathrm{sym}}$ with $\bm{0}\leq\bm{K}\leq\bm{K}^{\max}(\bm{L}+\bm{E}^{(ij)})$, and under the convention that $H_{\bm{L},\bm{K}}=0$ for $\bm{K}<\bm{0}$ or $\bm{K}>\bm{K}^{\max}(\bm{L})$, the recurrence relation for $H_{\bm{L},\bm{K}}$ is 
\begin{equation}\label{H:recurrence}
H_{\bm{L}+\bm{E}^{(ij)},\bm{K}}=H_{\bm{L},\bm{K}}+[\ell_{ji}-(1+\delta_{ij})(k_{ij}-1)]H_{\bm{L},\bm{K}-\bm{E}^{(ij)}_{\text{sym}}}.
\end{equation}
\end{lemma}
\begin{proof}
See Appendix \ref{A:H}.
\end{proof}
Thus, $A_j+B_j$ reads
\begin{align}\label{eq:A+B}
A_j+B_j=\sum_{\bm{r}(\bm{L})\leq\bm{n}} \binom{\bm{n}}{\bm{L}}_r\bm{\mu}^{\bm{n}-\bm{r}(\bm{L})}&\sum_{\bm{K}\leq\bm{K}^{\max}(\bm{L}+\bm{E}^{(ij)})}H_{\bm{L}+\bm{E}^{(ij)},\bm{K}}\bm{C}_U^{\tilde{\bm{L}}-\bm{K}+\bm{E}^{(ij)}_{\text{sym}}}\times\nonumber\\&\times\mathsf{E}\left[\partial_j^{c_j(\bm{L}-\tilde{\bm{K}})+1}\prod_{\substack{p=1\\ p\neq j}}^N\partial_p^{c_{p}(\bm{L}-\tilde{\bm{K}})}g(\bm{X})\right],
\end{align}
and by performing the index-matrix change $\bm{L}'=\bm{L}+\bm{E}^{(ij)}$
\begin{align}\label{eq:A+B2}
A_j+B_j=\sum_{\bm{r}(\bm{L}')\leq\bm{n}+\bm{e}^{(i)}} \binom{\bm{n}}{\bm{L}'-\bm{E}^{(ij)}}_r\bm{\mu}^{\bm{n}+\bm{e}^{(i)}-\bm{r}(\bm{L}')}&\sum_{\bm{K}\leq\bm{K}^{\max}(\bm{L}')}H_{\bm{L}',\bm{K}}\bm{C}_U^{\tilde{\bm{L}}'-\bm{K}}\times\nonumber\\&\times\mathsf{E}\left[\partial^{\bm{c}(\bm{L}'-\tilde{\bm{K}})}g(\bm{X})\right].
\end{align}
By also considering $A_0$ from Eq.~\eqref{eq:A_0}, Eq.~\eqref{eq:AB2} reads
\begin{align}\label{eq:last}
\mathsf{E}\left[(g(\bm{X})X_i)\bm{X}^{\bm{n}}\right]=&\sum_{\bm{r}(\bm{L})\leq\bm{n}} \left[\binom{\bm{n}}{\bm{L}}_r+\sum_{j=1}^N\binom{\bm{n}}{\bm{L}-\bm{E}^{(ij)}}_r\right]\bm{\mu}^{\bm{n}+\bm{e}^{(i)}-\bm{r}(\bm{L})}\times\nonumber\\&\times\sum_{\bm{K}\leq\bm{K}^{\max}(\bm{L})}H_{\bm{L},\bm{K}}\bm{C}_U^{\tilde{\bm{L}}-\bm{K}}\mathsf{E}\left[\partial^{\bm{c}(\bm{L}-\tilde{\bm{K}})}g(\bm{X})\right]+\nonumber\\
&+\sum_{\bm{r}(\bm{L})=\bm{n}+\bm{e}^{(i)}} \left[\sum_{j=1}^N\binom{\bm{n}}{\bm{L}-\bm{E}^{(ij)}}_r\right]\bm{\mu}^{\bm{n}+\bm{e}^{(i)}-\bm{r}(\bm{L})}\times\nonumber\\&\times\sum_{\bm{K}\leq\bm{K}^{\max}(\bm{L})}H_{\bm{L},\bm{K}}\bm{C}_U^{\tilde{\bm{L}}-\bm{K}}\mathsf{E}\left[\partial^{\bm{c}(\bm{L}-\tilde{\bm{K}})}g(\bm{X})\right].
\end{align}
The inductive proof of Eq.~\eqref{eq:stein_new} is completed by the following Lemma.
\begin{lemma}\label{lem:multi}
The addition of $r$-multinomial coefficients reads
\begin{equation}\label{eq:multi1}
\binom{\bm{n}+\bm{e}^{(i)}}{\bm{L}}_r=\binom{\bm{n}}{\bm{L}}_r+\sum_{j=1}^N\binom{\bm{n}}{\bm{L}-\bm{E}^{(ij)}}_r, \ \ \text{for} \ \ \bm{r}(\bm{L})\leq\bm{n},
\end{equation}
\begin{equation}\label{eq:multi2}
\binom{\bm{n}+\bm{e}^{(i)}}{\bm{L}}_r=\sum_{j=1}^N\binom{\bm{n}}{\bm{L}-\bm{E}^{(ij)}}_r, \ \ \text{for} \ \ \bm{r}(\bm{L})=\bm{n}+\bm{e}^{(i)}.
\end{equation}
\end{lemma}
\begin{proof}
See Appendix \ref{A:multi}.
\end{proof}
Substitution of Eqs.~\eqref{eq:multi1}, \eqref{eq:multi2} into Eq.~\eqref{eq:last} results Eq.~\eqref{eq:stein_new} for $\bm{n}+\bm{e}^{(i)}$.

\section{Constructive formal derivation of Theorem \ref{cor:index_matrix}}\label{sec:construction} 
Our alternative, constructive derivation of Eq.~\eqref{eq:stein_new} is based on the following relation for the mean value operator.
\begin{theorem}\label{def1}
Let $\bm{X}$ be an $N$-dimensional Gaussian random vector with mean value vector $\bm{\mu}$ and autocovariance matrix $\bm{C}$, and $g$ be a $C^{\infty}\left(\mathbb{R}^N\rightarrow\mathbb{R}\right)$ function. The diagonal elements of matrix $\bm{C}$ (the autocovariances of each $X_i$ component) are denoted as $\sigma^2_i=C_{ii}$. The expectation $\mathsf{E}\left[g(\bm{X})\right]$ is expressed as the action of averaged shift operators $\mathcal{T}_{ij}$:
\begin{equation}\label{eq:mean_value_action}
\mathsf{E}\left[g(\bm{X})\right]=\left(\prod_{i=1}^N\mathcal{T}_{ii}\right)\left(\prod_{i=1}^N\prod_{j>i}^N\mathcal{T}_{ij}\right)g(\bm{\mu}).
\end{equation} 
$\mathcal{T}_{ij}$ are the pseudodifferential operators
\begin{equation}\label{Tii}
\mathcal{T}_{ii}=\exp\left(\frac{\sigma^2_{i}}{2}\partial^2_i\right), \ \ i=1,\ldots,N,
\end{equation}
\begin{equation}\label{Tij}
\mathcal{T}_{ij}=\exp\left(C_{ij}\partial_i\partial_j\right), \ \ i,j=1,\ldots,N, \ j\neq i,
\end{equation}
whose action is to be understood by their series forms
\begin{equation}\label{Tii_series}
\mathcal{T}_{ii}=\sum_{m=0}^{\infty}\frac{\sigma^{2m}_{i}}{2^mm!}\partial_i^{2m}, \ \ i=1,\ldots,N,
\end{equation}
\begin{equation}\label{Tij_series}
\mathcal{T}_{ij}=\sum_{m=0}^{\infty}\frac{C_{ij}^m}{m!}\partial_i^m\partial_j^m, \ \ i,j=1,\ldots,N, \ j\neq i.
\end{equation}
\end{theorem}
\begin{proof}
We formally derive Eq.~\eqref{eq:mean_value_action} in Appendix \ref{A:def1}, using the moment-generating function of the $N$-dimensional Gaussian vector $\bm{X}$.  The infinite-dimensional counterpart of Theorem \ref{def1}, regarding Gaussian processes, is presented in \cite{mamis_NF}, and it is also found in \cite[Ch. 4]{klya} as a concept.
\end{proof}
\begin{remark}\label{rem:Tproperties}
Under the formal assumption that all infinite series involved are summable,  and by employing the linearity of derivatives, we can easily see that $\mathcal{T}_{ij}$ operators are linear, commute with differentiation operators $\partial_i$, and also commute with each other (see also \cite[Lemmata 1-3]{mamis_NF}).
\end{remark}
\begin{lemma}\label{lem:Tii} We determine the action of $\mathcal{T}_{ii}$ operator to
\begin{equation}\label{eq:Tii_action}
\mathcal{T}_{ii}\left[g(\bm{x})x_i^{n_i}\right]=\sum_{\ell=0}^{n_i}\binom{n_i}{\ell}x_i^{n_i-\ell}\sum_{k=0}^{\lfloor\ell/2\rfloor}H_{\ell,k}\sigma_i^{2(\ell-k)}\mathcal{T}_{ii}\left[\partial_i^{\ell-2k}g(\bm{x})\right],
\end{equation}
where $H_{\ell,k}$ are the signless Hermite coefficients, defined by Eq.~\eqref{eq:hermite_numbers}.
\end{lemma}
\begin{proof}
See Appendix \ref{A:Tii}.
\end{proof}
\begin{lemma}\label{lem:Tij} We determine the action of $\mathcal{T}_{ij}$, $j\neq i$ operator to
\begin{align}\label{eq:Tij_action}
\mathcal{T}_{ij}\left[g(\bm{x})x_i^{n_i}x_j^{n_j}\right]=\sum_{\ell_i=0}^{n_i}\binom{n_i}{\ell_i}x_i^{n_i-\ell_i}\sum_{\ell_j=0}^{n_j}\binom{n_j}{\ell_j}x_j^{n_j-\ell_j}&\sum_{k=0}^{\min\{\ell_i,\ell_j\}}G_{\ell_i,\ell_j,k}C_{ij}^{\ell_i+\ell_j-k}\times\nonumber\\&\times\mathcal{T}_{ij}\left[\partial_i^{\ell_j-k}\partial_j^{\ell_i-k}g(\bm{x})\right],
\end{align}
with
\begin{equation}\label{eq:glue_numbers}
G_{\ell_1,\ell_2,k}=\binom{\ell_1}{k}\binom{\ell_2}{k}k!, \ \ k=0,\ldots,\min\{\ell_1,\ell_2\}.
\end{equation}
Similarly to $H_{\ell,k}$ coefficients, $G_{\ell_1,\ell_2,k}$'s are identified as the absolute values of the coefficients appearing in the two-dimensional It{\^o}--Hermite polynomials \cite{ito,ismail}.
\end{lemma}
\begin{proof}
See Appendix \ref{A:Tij}.
\end{proof}
By expressing the expectation $\mathsf{E}\left[g\left(\bm{X}\right)\left(\prod_{i=1}^NX_i^{n_i}\right)\right]$ via Theorem \ref{def1},  we understand that, for its evaluation,  it suffices to sequentially apply operators $\mathcal{T}_{ii}$, $\mathcal{T}_{ij}$, $i,j=1,\ldots,N$, $j>i$ at the product $g\left(\bm{x}\right)\left(\prod_{i=1}^Nx_i^{n_i}\right)$, and set $\bm{x}=\bm{\mu}$ afterwards.  After algebraic manipulations and using the operator properties of remark \ref{rem:Tproperties}, we obtain 
\begin{align}\label{eq:formula}
&\mathsf{E}\left[g\left(\bm{X}\right)\left(\prod_{i=1}^NX_i^{n_i}\right)\right]=\sum_{\substack{m_i+\sum_{j=1}^N\ell_{ij}=n_i\\ i=1,\ldots,N}}\left[\prod_{i=1}^N\binom{n_i}{m_i,\ell_{i1},\ldots,\ell_{iN}}\right]\left(\prod_{i=1}^N\mu_i^{m_i}\right)\times\nonumber\\&\times\left(\prod_{i=1}^N\sum_{k_{ii}=0}^{\lfloor\ell_{ii}/2\rfloor}H_{\ell_{ii},k_{ii}}\sigma_{i}^{2(\ell_{ii}-k_{ii})}\right)\left(\prod_{i=1}^N\prod_{j>i}^N\sum_{k_{ij}=0}^{\min\{\ell_{ij},\ell_{ji}\}}G_{\ell_{ij},\ell_{ji},k_{ij}}C_{ij}^{\ell_{ij}+\ell_{ji}-k_{ij}}\right)\times\nonumber\\
&\times\mathsf{E}\left[\prod_{i=1}^N\partial_i^{\sum_{j=1}^N\left[\ell_{ji}-(1+\delta_{ij})k_{ij}\right]}g\left(\bm{X}\right)\right], \ \ \text{with} \ \  k_{ij}=k_{ji},
\end{align} 
where  $\binom{n}{m,\ell_{1},\ldots,\ell_{N}}=\frac{n!}{m!\ell_1!\cdots\ell_N!}$ is the multinomial coefficient with $N+1$ factors.
Sum $\sum_{\substack{m_i+\sum_{j=1}^N\ell_{ij}=n_i\\ i=1,\ldots,N}}$ is over all combinations of integers $\left\{\{m_i,\ell_{i1},\ldots,\ell_{iN}\}_{i=1}^N\right\}$ with $m_i+\sum_{j=1}^N\ell_{ij}=n_i$, $i=1,\ldots,N$. By recasting Eq.~\eqref{eq:formula} into multi-index and index matrix notation, we obtain Eq.~\eqref{eq:stein_new}.

\section{Conclusions and future work}\label{sec:conclusions}
In the present work, we derived formula \eqref{eq:stein_new} extending Stein's lemma for the evaluation of $\mathsf{E}\left[g(\bm{X})\prod_{i=1}^NX_i^{n_i}\right]$, where $\bm{X}$ is an $N$-dimensional Gaussian random vector.  By our formula,  the said expectation is expressed in terms of the expectations of partial derivatives of $g(\bm{X})$, as well as the mean value vector and autocovariance matrix of $\bm{X}$.  Furthermore, by setting $g(\bm{X})=1$, formula \eqref{eq:stein_new} results in Isserlis theorem \cite{isserlis} and Song \& Lee formula \cite{song} for Gaussian product moments $\mathsf{E}\left[\prod_{i=1}^NX_i^{n_i}\right]$.

A direction for future works is the extension of the infinite-dimensional analog of Stein's lemma, called the Novikov-Furutsu theorem (see \cite[Sec. 11.5]{scott}, \cite{mamis_NF}). In the infinite-dimensional case, $X$ is a Gaussian random process of time argument $t$, whose mean value is the function $\mu(t)$, and its two-time autocovariance function is $C(t_1,t_2)$.  Thus, for $g$ being a functional of $X$ over the time interval $[t_0,t]$, Novikov-Furutsu theorem reads:
\begin{equation}\label{eq:NF}
\mathsf{E}\left[g[X]X(t)\right]=\mu(t)\mathsf{E}\left[g(X)\right]+\int_{t_0}^tC(t,s)\mathsf{E}\left[\frac{\delta g[X]}{\delta X(s)}\right]\mathrm{d}s,
\end{equation}
where $\delta g[X]/\delta X(s)$ is the Volterra functional derivative of $g[X]$ with respect to a local perturbation of process $X$ centered at time $s$ (see e.g. \cite[Appendix A]{mamis_NF} for more on Volterra calculus).  Novikov-Furutsu theorem, Eq.~\eqref{eq:NF}, is the main tool in deriving evolution equations, that resemble the classical Fokker-Planck equation, for the response probability density of dynamical systems under Gaussian random excitation, see e.g. \cite[Eq.(3.19)]{hanggi}, \cite{Mamis2019,Mamis2021,Mamis2024}.  Recently \cite[Ch. 3]{Mamis2020}, we extended Novikov-Furutsu theorem for expectations that contain the Gaussian argument at various times; $\mathsf{E}\left[g[X]\prod_{i=1}^NX(t_i)\right]$. As we have already shown in \cite{mamis_NF}, the introduction and use of averaged shift operators is very helpful in constructing extensions of the Novikov-Furutsu theorem.

\appendix 
\section{Proof of Lemma \ref{lem:H}}\label{A:H}
By using the definition relation \eqref{eq:D}, we easily calculate that
\begin{equation}\label{eq:rec_0}
H_{\bm{L}+\bm{E}^{(ij)},\bm{0}}=H_{\bm{L},\bm{0}}=1.
\end{equation}
Since, by convention, $H_{\bm{L},-\bm{E}^{(ij)}_{\text{sym}}}=0$, Eq.~\eqref{eq:rec_0} coincides with recurrence relation \eqref{H:recurrence} for $\bm{K}=\bm{0}$. For $\bm{E}^{(ij)}_{\text{sym}}\leq \bm{K}\leq \bm{K}^{\max}(\bm{L})$, we have the following two cases.\\
\textit{First case:} $i=j$
\begin{align*}
&H_{\bm{L},\bm{K}}+(\ell_{ii}-2k_{ii}+2)H_{\bm{L},\bm{K}-\bm{E}^{(ii)}_{\text{sym}}}=\nonumber\\&=\frac{\bm{L}!}{2^{\sum_{p=1,p\neq i}^Nk_{pp}}\prod_{\substack{p=1,q=p\\(p,q)\neq(i,i)}}^Nk_{pq}!\prod_{\substack{p,q=1\\(p,q)\neq(i,i)}}^N(\ell_{pq}-\tilde{k}_{pq})!}\times\nonumber\\&\times\left[\frac{1}{2^{k_{ii}}k_{ii}!(\ell_{ii}-2k_{ii})!}+\frac{\ell_{ii}-2k_{ii}+2}{2^{k_{ii}-1}(k_{ii}-1)!(\ell_{ii}-2k_{ii}+2)!}\right]=\nonumber\\
&=\frac{\bm{L}!}{2^{\sum_{p=1,p\neq i}^Nk_{pp}}\prod_{\substack{p=1,q=p\\(p,q)\neq(i,i)}}^Nk_{pq}!\prod_{\substack{p,q=1\\(p,q)\neq(i,i)}}^N(\ell_{pq}-\tilde{k}_{pq})!}\frac{\ell_{ii}+1}{2^{k_{ii}}k_{ii}!(\ell_{ii}+1-2k_{ii})!}=\nonumber\\&=\frac{(\bm{L}+\bm{E}^{(ii)})!}{2^{\mathrm{tr}(\bm{K})}\bm{K}_U!(\bm{L}+\bm{E}^{(ii)}-\tilde{\bm{K}})!}=H_{\bm{L}+\bm{E}^{(ii)},\bm{K}}.
\end{align*}
\textit{Second case:} $i\neq j$
\begin{align*}
&H_{\bm{L},\bm{K}}+(\ell_{ji}-k_{ij}+1)H_{\bm{L},\bm{K}-\bm{E}^{(ij)}_{\text{sym}}}=\frac{\bm{L}!}{2^{\mathrm{tr}(\bm{K})}\prod_{\substack{p=1,q=p\\(p,q)\neq(i,j)}}^Nk_{pq}!\prod_{\substack{p,q=1\ \ \ \ \ \ \ \ \\(p,q)\neq(i,j),(j,i)}}^N(\ell_{pq}-\tilde{k}_{pq})!}\times\nonumber\\
&\times\left[\frac{1}{k_{ij}!(\ell_{ij}-k_{ij})!(\ell_{ji}-k_{ji})!}+\frac{\ell_{ji}-k_{ij}+1}{(k_{ij}-1)!(\ell_{ij}-k_{ij}+1)!(\ell_{ji}-k_{ji}+1)!}\right],
\end{align*}
and since $k_{ij}=k_{ji}$
\begin{align*}
&H_{\bm{L},\bm{K}}+(\ell_{ji}-k_{ij}+1)H_{\bm{L},\bm{K}-\bm{E}^{(ij)}_{\text{sym}}}=\frac{\bm{L}!}{2^{\mathrm{tr}(\bm{K})}\prod_{\substack{p=1,q=p\\(p,q)\neq(i,j)}}^Nk_{pq}!\prod_{\substack{p,q=1\ \ \ \ \ \ \ \ \\(p,q)\neq(i,j),(j,i)}}^N(\ell_{pq}-\tilde{k}_{pq})!}\times\nonumber\\&\times\frac{\ell_{ij}+1}{k_{ij}!(\ell_{ij}+1-k_{ij})!(\ell_{ji}-k_{ji})!}=\frac{(\bm{L}+\bm{E}^{(ij)})!}{2^{\mathrm{tr}(\bm{K})}\bm{K}_U!(\bm{L}+\bm{E}^{(ij)}-\tilde{\bm{K}})!}=H_{\bm{L}+\bm{E}^{(ij)},\bm{K}}.
\end{align*}
Last, we have to prove Eq.~\eqref{H:recurrence} for $\bm{K}=\bm{K}^{\max}(\bm{L}+\bm{E}^{(ij)})$. Again, we distinguish two cases.\\
\textit{First case:} $i=j$. If $\ell_{ii}$ is even, $\lfloor(\ell_{ii}+1)/2\rfloor=\lfloor\ell_{ii}/2\rfloor$, and thus $\bm{K}^{\max}(\bm{L}+\bm{E}^{(ii)})=\bm{K}^{\max}(\bm{L})$. So, it remains to prove Eq.~\eqref{H:recurrence} for odd $\ell_{ii}=2a+1$ and for $k_{ii}=k_{ii}^{\max}(\bm{L}+\bm{E}^{(ii)})=\lfloor(\ell_{ii}+1)\rfloor=a+1$. In this case, we calculate
\begin{align}\label{eq:rec_1}
\frac{H_{\bm{L}+\bm{E}^{(ii),\bm{K}^{\max}(\bm{L}+\bm{E}^{(ii)})}}}{H_{\bm{L},\bm{K}^{\max}(\bm{L}+\bm{E}^{(ii)})-\bm{E}^{(ii)}_{\text{sym}}}}=\frac{\frac{(2a+2)!}{2^{a+1}(a+1)!}}{\frac{(2a+1)!}{2^aa!}}=1.
\end{align}
Since, by convention, $H_{\bm{L},\bm{K}^{\max}(\bm{L}+\bm{E}^{(ii)})}=0$ for this case, Eq.~\eqref{eq:rec_1} is the specification of recurrence relation \eqref{H:recurrence}.\\
\textit{Second case:} $i\neq j$. If $\ell_{ij}+1\geq\ell_{ji}$, $\bm{K}^{\max}(\bm{L}+\bm{E}^{(ij)})=\bm{K}^{\max}(\bm{L})$. Thus, it remains to prove Eq.~\eqref{H:recurrence} for $\ell_{ij}+1<\ell_{ji}$, and for $k_{ij}=k_{ij}^{\max}(\bm{L}+\bm{E}^{(ij)})=\ell_{ij}+1$. In this case, we calculate
\begin{align}\label{eq:rec_2}
\frac{H_{\bm{L}+\bm{E}^{(ij),\bm{K}^{\max}(\bm{L}+\bm{E}^{(ij)})}}}{H_{\bm{L},\bm{K}^{\max}(\bm{L}+\bm{E}^{(ij)})-\bm{E}^{(ij)}_{\text{sym}}}}=\frac{\frac{(\ell_{ij}+1)!}{(\ell_{ij}+1)!(\ell_{ij}+1-\ell_{ij}-1)!(\ell_{ji}-\ell_{ij}-1)!}}{\frac{\ell_{ij}!}{(\ell_{ij}+1-1)!(\ell_{ij}-\ell_{ij}+1-1)!(\ell_{ji}-\ell_{ij})!}}=\ell_{ji}-\ell_{ij}.
\end{align}
Since, by convention, $H_{\bm{L},\bm{K}^{\max}(\bm{L}+\bm{E}^{(ij)})}=0$ for this case, Eq.~\eqref{eq:rec_2} is the specification of recurrence relation \eqref{H:recurrence}.

\section{Proof of Lemma \ref{lem:multi}}\label{A:multi}
\textit{First case:} $\bm{r}(\bm{L})\leq\bm{n}$. Using the definition relation \eqref{eq:multi_r} for $r$-multinomial coefficients, we have
\begin{align*}
\binom{\bm{n}}{\bm{L}}_r+\sum_{j=1}^N\binom{\bm{n}}{\bm{L}-\bm{E}^{(ij)}}_r&=\frac{\bm{n}!}{[\bm{n}-\bm{r}(\bm{L})]!\bm{L}!}+\sum_{j=1}^N\frac{\bm{n}!}{[\bm{n}-\bm{r}(\bm{L}-\bm{E}^{(ij)})]![\bm{L}-\bm{E}^{(ij)}]!}=\nonumber\\&=\frac{\bm{n}!}{[\bm{n}-\bm{r}(\bm{L})]!\bm{L}!}+\sum_{j=1}^N\frac{\bm{n}!}{[\bm{n}+\bm{e}^{(i)}-\bm{r}(\bm{L})]![\bm{L}-\bm{E}^{(ij)}]!}=\nonumber\\&=\frac{\bm{n}!}{[\bm{n}+\bm{e}^{(i)}-\bm{r}(\bm{L})]!\bm{L}!}\left[n_i+1-r_i(\bm{L})+\sum_{j=1}^N\ell_{ij}\right]=\nonumber\\&=\frac{[\bm{n}+\bm{e}^{(i)}]!}{[\bm{n}+\bm{e}^{(i)}-\bm{r}(\bm{L})]!\bm{L}!}=\binom{\bm{n}+\bm{e}^{(i)}}{\bm{L}}_r.
\end{align*}
\textit{Second case:} $\bm{r}(\bm{L})=\bm{n}+\bm{e}^{(i)}$. Using the definition relation \eqref{eq:multi_r}, we have
\begin{align*}
\sum_{j=1}^N\binom{\bm{n}}{\bm{L}-\bm{E}^{(ij)}}_r&=\sum_{j=1}^N\frac{\bm{n}!}{[\bm{n}-\bm{r}(\bm{L}-\bm{E}^{(ij)})]![\bm{L}-\bm{E}^{(ij)}]!}=\nonumber\\&=\sum_{j=1}^N\frac{\bm{n}!}{[\bm{n}+\bm{e}^{(i)}-\bm{r}(\bm{L})]![\bm{L}-\bm{E}^{(ij)}]!}=\frac{\bm{n}!\sum_{j=1}^N\ell_{ij}}{[\bm{n}+\bm{e}^{(i)}-\bm{r}(\bm{L})]!\bm{L}!}=\nonumber\\&=\frac{\bm{n}!(n_i+1)}{[\bm{n}+\bm{e}^{(i)}-\bm{r}(\bm{L})]!\bm{L}!}=\frac{[\bm{n}+\bm{e}^{(i)}]!}{[\bm{n}+\bm{e}^{(i)}-\bm{r}(\bm{L})]!\bm{L}!}=\binom{\bm{n}+\bm{e}^{(i)}}{\bm{L}}_r.
\end{align*}

\section{Formal derivation of Theorem \ref{def1}}\label{A:def1}
The Taylor expansion of a $C^{\infty}\left(\mathbb{R}^N\rightarrow\mathbb{R}\right)$ function $g$ around $\bm{x}_0$ is expressed via the shift pseudodifferential operator in exponential form (see e.g. \cite[Sec. 1.1]{glaeske}) as
\begin{align}\label{eq:shift}
g(\bm{x})&=\left(1+\sum_{m=1}^{\infty}\frac{1}{m!}\sum_{i_1=1}^N\stackrel{(m)}{\cdots}\sum_{i_m=1}^N\hat{x}_{i_1}\cdots\hat{x}_{i_m}\partial_{i_1}\cdots\partial_{i_m}\right)g(\bm{x}_0)=\nonumber\\&=\left[\sum_{m=0}^{\infty}\frac{1}{m!}\left(\sum_{i=1}^N\hat{x}_i\partial_i\right)^m\right]g(\bm{x}_0)=\exp\left(\sum_{i=1}^N\hat{x}_i\partial_i\right)g(\bm{x}_0),
\end{align}
where $\hat{\bm{x}}=\bm{x}-\bm{x}_0$ is called the shift argument. By substituting the random vector $\bm{X}$ as the argument of function $g$, choosing $\bm{x}_0=\bm{\mu}$ where $\bm{\mu}$ is the mean value of $\bm{X}$, and taking the expectation in both sides of Eq.~\eqref{eq:shift} results into
\begin{equation}\label{eq:averaged_shift}
\mathsf{E}[g(\bm{X})]=\mathsf{E}\left[\exp\left(\sum_{i=1}^N\hat{X}_i\partial_i\right)\right]g(\bm{\mu})=M_{\bm{\hat{X}}}(\bm{\nabla})g(\bm{\mu}),
\end{equation}
where $\bm{\hat{X}}:=\bm{X}-\bm{\mu}$ (the centered random vector) and $\bm{\nabla}=[\partial_1,\ldots,\partial_N]^T$ (the del vector).  In Eq.~\eqref{eq:averaged_shift}, $M_{\bm{\hat{X}}}(\bm{u})$ is identified as the moment-generating function of $\bm{\hat{X}}$; $M_{\bm{\hat{X}}}(\bm{u})=\mathsf{E}\left[\exp\left(\bm{u}^T\bm{\hat{X}}\right)\right]$, see \cite[Sec. 4.3.3]{song2}. For the Gaussian vector $\bm{X}$ with autocovariance matrix $\bm{C}$, the moment-generating function for the corresponding centered Gaussian random vector $\bm{\hat{X}}$ takes the form $M_{\bm{\hat{X}}}(\bm{u})=\exp\left(\bm{u}^T\bm{C}\bm{u}/2\right)$, see \cite[Sec. 5.1.1]{song2}.  Substitution of Gaussian $M_{\bm{\hat{X}}}(\bm{u})$ into Eq.~\eqref{eq:averaged_shift} results in
\begin{align*}
\mathsf{E}[g(\bm{X})]=\exp\left(\frac{1}{2}\sum_{i=1}^N\sum_{j=1}^NC_{ij}\partial_i\partial_j\right)g(\bm{\mu}),
\end{align*}
and by using the symmetry property of autocovariance matrix $\bm{C}$:
\begin{align}\label{eq:gaussian_averaged_shift1}
\mathsf{E}[g(\bm{X})]&=\exp\left(\sum_{i=1}^N\frac{\sigma^2_i}{2}\partial^2_i+\sum_{i=1}^N\sum_{j>i}^NC_{ij}\partial_i\partial_j\right)g(\bm{\mu})=\nonumber\\&=\left[\prod_{i=1}^N\exp\left(\frac{\sigma^2_i}{2}\partial^2_i\right)\right]\left[\prod_{i=1}^N\prod_{j>i}^N\exp\left(C_{ij}\partial_i\partial_j\right)\right]g(\bm{\mu}).
\end{align}
Eq.~\eqref{eq:gaussian_averaged_shift1} coincides with Eq.~\eqref{eq:mean_value_action}.

\section{Proof of Lemma \ref{lem:Tii}}\label{A:Tii}
Expressing $\mathcal{T}_{ii}\left[g(\bm{x})x_i^{n_i}\right]$ via Eq.~\eqref{Tii_series} we have
\begin{equation}\label{eq:Tii1}
\mathcal{T}_{ii}\left[g(\bm{x})x_i^{n_i}\right]=\sum_{m=0}^{\infty}\frac{\sigma^{2m}_{i}}{2^mm!}\partial_i^{2m}\left(g(\bm{x})x_i^{n_i}\right).
\end{equation}
The derivatives appearing in the right-hand side of Eq.~\eqref{eq:Tii1} are further evaluated using the general Leibniz rule:
\begin{equation}\label{gen_leibniz}
\partial_i^{2m}\left(g(\bm{x})x_i^{n_i}\right)=\sum_{\ell=0}^{2m}\binom{2m}{\ell}\left(\partial_i^{2m-\ell}g(\bm{x})\right)\left(\partial_i^{\ell}x_i^{n_i}\right).
\end{equation}
Since $\partial_i^{\ell}x_i^{n_i}=(n_i!/(n_i-\ell)!)x_i^{n_i-\ell}$ for $n_i\geq\ell$ and zero for $n_i<\ell$, Eq.~\eqref{gen_leibniz} is updated to
\begin{equation}\label{gen_leibniz1}
\partial_i^{2m}\left(g(\bm{x})x_i^{n_i}\right)=\sum_{\ell=0}^{\min\{n_i,2m\}}\binom{n_i}{\ell}(2m)^{\underline{\ell}}\left(\partial_i^{2m-\ell}g(\bm{x})\right)x_i^{n_i-\ell},
\end{equation}
where $(2m)^{\underline{\ell}}=(2m)(2m-1)\cdots(2m-\ell+1)$ is the falling factorial.  Substitution of Eq.~\eqref{gen_leibniz1} into Eq.~\eqref{eq:Tii1} results in
\begin{equation}\label{eq:Tii2}
\mathcal{T}_{ii}\left[g(\bm{x})x_i^{n_i}\right]=\sum_{m=0}^{\infty}\sum_{\ell=0}^{\min\{n_i,2m\}}\binom{n_i}{\ell}x_i^{n_i-\ell}\frac{\sigma^{2m}_{i}}{2^m}\frac{(2m)^{\underline{\ell}}}{m!}\partial_i^{2m-\ell}g(\bm{x}).
\end{equation}
In Eq.~\eqref{eq:Tii2}, the $m$ and $\ell$-summations are interchanged using formula \eqref{sum1}, resulting into
\begin{equation}\label{eq:Tii3}
\mathcal{T}_{ii}\left[g(\bm{x})x_i^{n_i}\right]=\sum_{\ell=0}^{n_i}\binom{n_i}{\ell}x_i^{n_i-\ell}\sum_{m=\lceil \ell/2\rceil}^{\infty}\frac{\sigma^{2m}_{i}}{2^m}\frac{(2m)^{\underline{\ell}}}{m!}\partial_i^{2m-\ell}g(\bm{x}).
\end{equation}
Following \cite[Sec. 8.4]{char}, see also \cite[Eq.(30)]{mamis_stein},  $(2m)^{\underline{\ell}}$ is expressed in terms of $m^{\underline{p}}$ as
\begin{equation*}
(2m)^{\underline{\ell}}=\sum_{p=\lceil\ell/2\rceil}^{\min\{m,\ell\}}C(\ell,p;2)m^{\underline{p}}
\end{equation*}
where $C(\ell,p;2)$ are the \emph{generalized factorial coefficients with parameter 2}, and $\lceil\bm{\cdot}\rceil$ is the ceiling function. Using also the fact $m^{\underline{p}}/m!=1/(m-p)!$, Eq.~\eqref{eq:Tii3} is expressed as
\begin{equation*}
\mathcal{T}_{ii}\left[g(\bm{x})x_i^{n_i}\right]=\sum_{\ell=0}^{n_i}\binom{n_i}{\ell}x_i^{n_i-\ell}\sum_{m=\lceil \ell/2\rceil}^{\infty}\sum_{p=\lceil\ell/2\rceil}^{\min\{m,\ell\}}\frac{\sigma^{2m}_{i}}{2^m}\frac{C(\ell,p;2)}{(m-p)!}\partial_i^{2m-\ell}g(\bm{x}).
\end{equation*}
By also interchanging the $m$ and $p$-summations using formula \eqref{sum2}, we have
\begin{equation*}
\mathcal{T}_{ii}\left[g(\bm{x})x_i^{n_i}\right]=\sum_{\ell=0}^{n_i}\binom{n_i}{\ell}x_i^{n_i-\ell}\sum_{p=\lceil \ell/2\rceil}^{\ell}\sum_{m=p}^{\infty}\frac{\sigma^{2m}_{i}}{2^m}\frac{C(\ell,p;2)}{(m-p)!}\partial_i^{2m-\ell}g(\bm{x}).
\end{equation*}
An index change in the $m$-sum, and the use of Eq.~\eqref{Tii_series}, results in
\begin{align*}
\mathcal{T}_{ii}\left[g(\bm{x})x_i^{n_i}\right]&=\sum_{\ell=0}^{n_i}\binom{n_i}{\ell}x_i^{n_i-\ell}\sum_{p=\lceil \ell/2\rceil}^{\ell}\frac{C(\ell,p;2)}{2^p}\sigma_i^{2p}\sum_{m=0}^{\infty}\frac{\sigma^{2m}_{i}}{2^mm!}\partial_i^{2m+2p-\ell}g(\bm{x})=\nonumber\\
&=\sum_{\ell=0}^{n_i}\binom{n_i}{\ell}x_i^{n_i-\ell}\sum_{p=\lceil \ell/2\rceil}^{\ell}\frac{C(\ell,p;2)}{2^p}\sigma_i^{2p}\mathcal{T}_{ii}\left[\partial_i^{2p-\ell}g(\bm{x})\right].
\end{align*}
Last, we perform the change of index $k=\ell-p$ to obtain
\begin{equation}\label{eq:Tii6}
\mathcal{T}_{ii}\left[g(\bm{x})x_i^{n_i}\right]=\sum_{\ell=0}^{n_i}\binom{n_i}{\ell}x_i^{n_i-\ell}\sum_{k=0}^{\lfloor\ell/2\rfloor}\frac{C(\ell,\ell-k;2)}{2^{\ell-k}}\sigma_i^{2(\ell-k)}\mathcal{T}_{ii}\left[\partial_i^{\ell-2k}g(\bm{x})\right].
\end{equation}
As we have showed in the recent work \cite[Lemma 2]{mamis_stein},  $H_{\ell,k}=C(\ell,\ell-k;2)/2^{\ell-k}$, for $\ell\in\mathbb{N}$, $k=0,\ldots,\lfloor\ell/2\rfloor$. Under this, Eq.~\eqref{eq:Tii6} coincides with Eq.~\eqref{eq:Tii_action}. For more on the connection between generalized factorial coefficients and Hermite polynomials, see \cite{char2}.

\section{Proof of Lemma \ref{lem:Tij}}\label{A:Tij}
Expressing $\mathcal{T}_{ij}\left[g(\bm{x})x_i^{n_i}x_j^{n_j}\right]$ via Eq.~\eqref{eq:Tij_action} we have
\begin{equation}\label{eq:Tij1}
\mathcal{T}_{ij}\left[g(\bm{x})x_i^{n_i}x_j^{n_j}\right]=\sum_{m=0}^{\infty}\frac{C_{ij}^m}{m!}\partial_i^m\partial_j^m\left(g(\bm{x})x_i^{n_i}x_j^{n_j}\right).
\end{equation}
As in Appendix \ref{A:Tii}, the derivatives $\partial_i^m$, $\partial_j^m$ in the right-hand side of Eq.~\eqref{eq:Tij1} can be evaluated further using general Leibniz rule, resulting in
\begin{align}\label{eq:Tij2}
\mathcal{T}_{ij}\left[g(\bm{x})x_i^{n_i}x_j^{n_j}\right]=&\sum_{m=0}^{\infty}\sum_{\ell_i=0}^{\min\{n_i,m\}}\binom{n_i}{\ell_i}x_i^{n_i-\ell_i}\times\nonumber\\&\times\sum_{\ell_j=0}^{\min\{n_j,m\}}\binom{n_j}{\ell_j}x_j^{n_j-\ell_j}C_{ij}^m\frac{m^{\underline{\ell_i}}m^{\underline{\ell_j}}}{m!}\partial_i^{m-\ell_i}\partial_j^{m-\ell_j}g(\bm{x}).
\end{align}
By rearranging the summations in Eq.~\eqref{eq:Tij2} using formula \eqref{sum3}, we obtain
\begin{align}\label{eq:Tij22}
\mathcal{T}_{ij}\left[g(\bm{x})x_i^{n_i}x_j^{n_j}\right]=&\sum_{\ell_i=0}^{n_i}\binom{n_i}{\ell_i}x_i^{n_i-\ell_i}\sum_{\ell_j=0}^{n_j}\binom{n_j}{\ell_j}x_j^{n_j-\ell_j}\times\nonumber\\&\times\sum_{m=\max\{\ell_i,\ell_j\}}^{\infty}C_{ij}^m\frac{m^{\underline{\ell_i}}m^{\underline{\ell_j}}}{m!}\partial_i^{m-\ell_i}\partial_j^{m-\ell_j}g(\bm{x}).
\end{align}
In order to evaluate the right-hand side of Eq.~\eqref{eq:Tij22} further, the product of the two falling factorials $m^{\underline{\ell_i}}m^{\underline{\ell_j}}$ has to be expressed in terms of falling factorials of $m$. This is performed by the following Lemma.
\begin{lemma}[Product of two falling factorials of $m$]\label{lem:falling}
It holds true that
\begin{equation}\label{G_property}
m^{\underline{\ell_i}}m^{\underline{\ell_j}}=\sum_{k=0}^{\min\{\ell_i,\ell_j\}}G_{\ell_i,\ell_j,k}m^{\underline{\ell_i+\ell_j-k}}.
\end{equation}
\end{lemma}
\begin{proof}
See Appendix \ref{A:falling}.
\end{proof}
Since, by the definition of falling factorial, $m^{\underline{\ell_i+\ell_j-k}}$ is zero for $\ell_i+\ell_j-k>m$,  Eq.~\eqref{G_property} is updated to
\begin{equation}\label{G_property2}
m^{\underline{\ell_i}}m^{\underline{\ell_j}}=\sum_{k=\max\{0,\ell_i+\ell_j-m\}}^{\min\{\ell_i,\ell_j\}}G_{\ell_i,\ell_j,k}m^{\underline{\ell_i+\ell_j-k}}.
\end{equation}
Substitution of Eq.~\eqref{G_property2} into Eq.~\eqref{eq:Tij22}, and use of $m^{\underline{\ell_i+\ell_j-k}}/m!=1/(m-\ell_i-\ell_j+k)!$ results in
\begin{align*}
\mathcal{T}_{ij}\left[g(\bm{x})x_i^{n_i}x_j^{n_j}\right]=&\sum_{\ell_i=0}^{n_i}\binom{n_i}{\ell_i}x_i^{n_i-\ell_i}\sum_{\ell_j=0}^{n_j}\binom{n_j}{\ell_j}x_j^{n_j-\ell_j}\times\nonumber\\&\times\sum_{m=\max\{\ell_i,\ell_j\}}^{\infty}\sum_{k=\max\{0,\ell_i+\ell_j-m\}}^{\min\{\ell_i,\ell_j\}}\frac{C_{ij}^mG_{\ell_i,\ell_j,k}}{(m-\ell_i-\ell_j+k)!}\partial_i^{m-\ell_i}\partial_j^{m-\ell_j}g(\bm{x}).
\end{align*}
By interchanging $m$ and $k$-summations using formula \eqref{sum4}, we have
\begin{align*}
\mathcal{T}_{ij}\left[g(\bm{x})x_i^{n_i}x_j^{n_j}\right]=&\sum_{\ell_i=0}^{n_i}\binom{n_i}{\ell_i}x_i^{n_i-\ell_i}\sum_{\ell_j=0}^{n_j}\binom{n_j}{\ell_j}x_j^{n_j-\ell_j}\times\nonumber\\&\times\sum_{k=0}^{\min\{\ell_i,\ell_j\}}G_{\ell_i,\ell_j,k}\sum_{m=\ell_i+\ell_j-k}^{\infty}\frac{C_{ij}^m}{(m-\ell_i-\ell_j+k)!}\partial_i^{m-\ell_i}\partial_j^{m-\ell_j}g(\bm{x}).
\end{align*}
An index change in the $m$-sum results in
\begin{align}\label{eq:Tij5}
\mathcal{T}_{ij}\left[g(\bm{x})x_i^{n_i}x_j^{n_j}\right]=&\sum_{\ell_i=0}^{n_i}\binom{n_i}{\ell_i}x_i^{n_i-\ell_i}\sum_{\ell_j=0}^{n_j}\binom{n_j}{\ell_j}x_j^{n_j-\ell_j}\times\nonumber\\&\times\sum_{k=0}^{\min\{\ell_i,\ell_j\}}G_{\ell_i,\ell_j,k}C_{ij}^{\ell_i+\ell_j-k}\sum_{m=0}^{\infty}\frac{C_{ij}^m}{m!}\partial_i^{m+\ell_j-k}\partial_j^{m+\ell_i-k}g(\bm{x}).
\end{align}
Eq.~\eqref{eq:Tij5}, by virtue of Eq.~\eqref{Tij_series}, coincides with Eq.~\eqref{eq:Tij_action}.

\section{Proof of Lemma \ref{lem:falling}}\label{A:falling}
Our starting point for proving Eq.~\eqref{G_property} is Vandermonde's identity (see e.g. \cite[example 3.6]{char})
\begin{equation}\label{eq:van}
\binom{m}{\ell_i}=\sum_{k=0}^{\ell_i}\binom{\ell_j}{k}\binom{m-\ell_j}{\ell_i-k}.
\end{equation}
By multiplying both sides of Eq.~\eqref{eq:van} by $\binom{m}{\ell_j}$, and after some algebraic manipulations, we obtain
\begin{equation}\label{eq:van2}
\binom{m}{\ell_i}\binom{m}{\ell_j}=\sum_{k=0}^{\ell_i}\binom{m}{\ell_j}\binom{\ell_j}{k}\binom{m-\ell_j}{\ell_i-k}=\sum_{k=0}^{\ell_i}\binom{\ell_i+\ell_j-k}{k,\ell_i-k,\ell_j-k}\binom{m}{\ell_i+\ell_j-k}.
\end{equation}
Since $\binom{\ell_i+\ell_j-k}{k,\ell_i-k,\ell_j-k}=0$ for $k>\ell_i$ or $k>\ell_j$, upper limit of $k$-sum in Eq.~\eqref{eq:van2} is updated to $\min\{\ell_i,\ell_j\}$.  By also using the fact that $m^{\underline{k}}=\binom{m}{k}k!$, see \cite[Eq. (3.11)]{char}, Eq.~\eqref{eq:van2} results in
\begin{equation}\label{eq:van3}
m^{\underline{\ell_i}}m^{\underline{\ell_j}}=\sum_{k=0}^{\min\{\ell_i,\ell_j\}}G_{\ell_i,\ell_j,k}m^{\underline{\ell_i+\ell_j-k}},
\end{equation}
where
\begin{equation}\label{eq:van4}
G_{\ell_i,\ell_j,k}=\binom{\ell_i+\ell_j-k}{k,\ell_i-k,\ell_j-k}\frac{\ell_i!\ell_j!}{(\ell_i+\ell_j-k)!}=\binom{\ell_i}{k}\binom{\ell_j}{k}k!.
\end{equation}
Eq.~\eqref{eq:van3} coincides with Eq.~\eqref{G_property}, and Eq.~\eqref{eq:van4} coincides with the definition relation \eqref{eq:glue_numbers} of $G_{\ell_i,\ell_j,k}$, completing thus the proof of Lemma \ref{lem:falling}.

\section{Summation rearrangement formulas and their proofs}
In this Appendix, we prove the formulas \eqref{sum1}, \eqref{sum2}, \eqref{sum3}, \eqref{sum4} rearranging multiple summations, that are employed in Appendices \ref{A:Tii}, \ref{A:Tij} for the proofs of Lemmata \ref{lem:Tii}, \ref{lem:Tij} respectively.
\begin{equation}\label{sum1}
\sum_{m=0}^{\infty}\sum_{\ell=0}^{\min\{n,2m\}}A_{m,\ell}=\sum_{\ell=0}^n\sum_{m=\lceil\ell/2\rceil}^{\infty}A_{m,\ell}.
\end{equation}
\begin{proof}
We distinguish the cases of even and odd $n$.  For $n=2p$, the left-hand side of Eq.~\eqref{sum1} is split into
\begin{equation}\label{sum1_1}
\sum_{m=0}^{\infty}\sum_{\ell=0}^{\min\{2p,2m\}}A_{m,\ell}=\sum_{m=0}^{p}\sum_{\ell=0}^{2m}A_{m,\ell}+\sum_{m=p+1}^{\infty}\sum_{\ell=0}^{2p}A_{m,\ell}.
\end{equation}
In the right-hand side of Eq.~\eqref{sum1_1}, the sums of the second term are interchanged. The double summation of the first term is over the triangle $0\leq m\leq p$, $0\leq \ell\leq 2m$ which is rearranged into $0\leq\ell\leq 2p$, $\ell/2\leq m\leq p$ and since $m$, $\ell$ are integers; $0\leq\ell\leq 2p$, $\lceil\ell/2\rceil\leq m\leq p$. Thus:
\begin{equation*}
\sum_{m=0}^{\infty}\sum_{\ell=0}^{\min\{2p,2m\}}A_{m,\ell}=\sum_{\ell=0}^{2p}\sum_{m=\lceil\ell/2\rceil}^{p}A_{m,\ell}+\sum_{\ell=0}^{2p}\sum_{m=p+1}^{\infty}A_{m,\ell}=\sum_{\ell=0}^{2p}\sum_{m=\lceil\ell/2\rceil}^{\infty}A_{m,\ell}
\end{equation*}
For $n=2p+1$, the left-hand side of Eq.~\eqref{sum1} is split into
\begin{align}\label{sum1_3}
\sum_{m=0}^{\infty}\sum_{\ell=0}^{\min\{2p+1,2m\}}A_{m,\ell}&=\sum_{m=0}^{p}\sum_{\ell=0}^{2m}A_{m,\ell}+\sum_{m=p+1}^{\infty}\sum_{\ell=0}^{2p+1}A_{m,\ell}=\nonumber\\&=\sum_{m=0}^{p}\sum_{\ell=0}^{2m}A_{m,\ell}+\sum_{m=p+1}^{\infty}\sum_{\ell=0}^{2p}A_{m,\ell}+\sum_{m=p+1}^{\infty}A_{m,2p+1}.
\end{align}
In the rightmost side of Eq.~\eqref{sum1_3}, the double summations are rearranged as for Eq.~\eqref{sum1_1}, resulting into
\begin{equation}\label{sum1_4}
\sum_{m=0}^{\infty}\sum_{\ell=0}^{\min\{2p+1,2m\}}A_{m,\ell}=\sum_{\ell=0}^{2p}\sum_{m=\lceil\ell/2\rceil}^{\infty}A_{m,\ell}+\sum_{m=p+1}^{\infty}A_{m,2p+1}.
\end{equation}
Since $p+1=\lceil(2p+1)/2\rceil$,  we identify the single sum in the right-hand side of Eq.~\eqref{sum1_4} as the $\ell=2p+1$ term of the double sum:
\begin{equation*}
\sum_{m=0}^{\infty}\sum_{\ell=0}^{\min\{2p+1,2m\}}A_{m,\ell}=\sum_{\ell=0}^{2p+1}\sum_{m=\lceil\ell/2\rceil}^{\infty}A_{m,\ell}.
\end{equation*}
Thus, the proof of Eq.~\eqref{sum1} for both even and odd $n$ is completed.  
\end{proof}
\begin{equation}\label{sum2}
\sum_{m=n}^{\infty}\sum_{k=n}^{\min\{\ell,m\}}A_{m,k}=\sum_{k=n}^{\ell}\sum_{m=k}^{\infty}A_{m,k}.
\end{equation}
\begin{proof}
The left-hand side of Eq.~\eqref{sum2} is split into
\begin{equation}\label{sum2_2}
\sum_{m=n}^{\infty}\sum_{k=n}^{\min\{\ell,m\}}A_{m,k}=\sum_{m=n}^{\ell}\sum_{k=n}^mA_{m,k}+\sum_{m=\ell+1}^{\infty}\sum_{k=n}^{\ell}A_{m,k}.
\end{equation}
In the right-hand side of Eq.~\eqref{sum2_2}, the sums of the second term are interchanged. The double summation of the first term is over the triangle $n\leq m\leq \ell$, $n\leq k\leq m$ which is rearranged into $n\leq k\leq \ell$, $k\leq m\leq \ell$. Thus:
\begin{equation*}
\sum_{m=n}^{\infty}\sum_{k=n}^{\min\{\ell,m\}}A_{m,k}=\sum_{k=n}^{\ell}\sum_{m=k}^{\ell}A_{m,k}+\sum_{k=n}^{\ell}\sum_{m=\ell+1}^{\infty}A_{m,k}=\sum_{k=n}^{\ell}\sum_{m=k}^{\infty}A_{m,k},
\end{equation*}
which completes the proof of Eq.~\eqref{sum2}.
\end{proof}
\begin{equation}\label{sum3}
\sum_{m=0}^{\infty}\sum_{\ell_1=0}^{\min\{n_1,m\}}\sum_{\ell_2=0}^{\min\{n_2,m\}}A_{m,\ell_1,\ell_2}=\sum_{\ell_1=0}^{n_1}\sum_{\ell_2=0}^{n_2}\sum_{m=\max\{\ell_1,\ell_2\}}^{\infty}A_{m,\ell_1,\ell_2}.
\end{equation}
\begin{proof}
Without loss of generality, we assume that $n_1<n_2$.  Then, the left-hand side of Eq.~\eqref{sum3} is split into
\begin{align}\label{sum3_2}
&\sum_{m=0}^{\infty}\sum_{\ell_1=0}^{\min\{n_1,m\}}\sum_{\ell_2=0}^{\min\{n_2,m\}}A_{m,\ell_1,\ell_2}=\nonumber\\&=\sum_{m=0}^{n_1}\sum_{\ell_1=0}^m\sum_{\ell_2=0}^mA_{m,\ell_1,\ell_2}+\sum_{m=n_1+1}^{n_2}\sum_{\ell_1=0}^{n_1}\sum_{\ell_2=0}^mA_{m,\ell_1,\ell_2}+\sum_{m=n_2+1}^{\infty}\sum_{\ell_1=0}^{n_1}\sum_{\ell_2=0}^{n_2}A_{m,\ell_1,\ell_2}.
\end{align}
In the first term of the right-hand side of Eq.~\eqref{sum3_2}, the summation is over $0\leq m\leq n_1$, $0\leq \ell_1\leq m$, $0\leq \ell_2\leq m$, which can be rearranged into $0\leq \ell_1\leq n_1$, $0\leq\ell_2\leq n_1$, $\max\{\ell_1,\ell_2\}\leq m\leq n_1$. In the second term, the summation over $n_1+1\leq m\leq n_2$, $0\leq\ell_2\leq m$ is rearranged into $0\leq\ell_2\leq n_2$, $\max\{n_1+1,\ell_2\}\leq m\leq n_2$.  Thus, Eq.~\eqref{sum3_2} is expressed equivalently as
\begin{align}\label{sum3_22}
&\sum_{m=0}^{\infty}\sum_{\ell_1=0}^{\min\{n_1,m\}}\sum_{\ell_2=0}^{\min\{n_2,m\}}A_{m,\ell_1,\ell_2}=\sum_{\ell_1=0}^{n_1}\sum_{\ell_2=0}^{n_1}\sum_{m=\max\{\ell_1,\ell_2\}}^{n_1}A_{m,\ell_1,\ell_2}+\nonumber\\&+\sum_{\ell_1=0}^{n_1}\sum_{\ell_2=0}^{n_2}\sum_{m=\max\{n_1+1,\ell_2\}}^{n_2}A_{m,\ell_1,\ell_2}+\sum_{\ell_1=0}^{n_1}\sum_{\ell_2=0}^{n_2}\sum_{m=n_2+1}^{\infty}A_{m,\ell_1,\ell_2}.
\end{align}
The second term in the right-hand side of Eq.~\eqref{sum3_22} is regrouped as
\begin{align}\label{sum3_3}
&\sum_{\ell_1=0}^{n_1}\sum_{\ell_2=0}^{n_2}\sum_{m=\max\{n_1+1,\ell_2\}}^{n_2}A_{m,\ell_1,\ell_2}=\nonumber\\&=\sum_{\ell_1=0}^{n_1}\sum_{\ell_2=0}^{n_1}\sum_{m=n_1+1}^{n_2}A_{m,\ell_1,\ell_2}+\sum_{\ell_1=0}^{n_1}\sum_{\ell_2=n_1+1}^{n_2}\sum_{m=\ell_2}^{n_2}A_{m,\ell_1,\ell_2}=\nonumber\\&=\sum_{\ell_1=0}^{n_1}\sum_{\ell_2=0}^{n_1}\sum_{m=n_1+1}^{n_2}A_{m,\ell_1,\ell_2}+\sum_{\ell_1=0}^{n_1}\sum_{\ell_2=n_1+1}^{n_2}\sum_{m=\ell_2}^{n_1}A_{m,\ell_1,\ell_2}+\sum_{\ell_1=0}^{n_1}\sum_{\ell_2=n_1+1}^{n_2}\sum_{m=n_1+1}^{n_2}A_{m,\ell_1,\ell_2}=\nonumber\\&=\sum_{\ell_1=0}^{n_1}\sum_{\ell_2=0}^{n_2}\sum_{m=n_1+1}^{n_2}A_{m,\ell_1,\ell_2}+\sum_{\ell_1=0}^{n_1}\sum_{\ell_2=n_1+1}^{n_2}\sum_{m=\ell_2}^{n_1}A_{m,\ell_1,\ell_2}=\nonumber\\
&=\sum_{\ell_1=0}^{n_1}\sum_{\ell_2=0}^{n_2}\sum_{m=n_1+1}^{n_2}A_{m,\ell_1,\ell_2}+\sum_{\ell_1=0}^{n_1}\sum_{\ell_2=n_1+1}^{n_2}\sum_{m=\max\{\ell_1\ell_2\}}^{n_1}A_{m,\ell_1,\ell_2}.
\end{align}
Substitution of Eq.~\eqref{sum3_3} into Eq.~\eqref{sum3_22} results into
\begin{align*}
&\sum_{m=0}^{\infty}\sum_{\ell_1=0}^{\min\{n_1,m\}}\sum_{\ell_2=0}^{\min\{n_2,m\}}A_{m,\ell_1,\ell_2}=\nonumber\\&=\sum_{\ell_1=0}^{n_1}\sum_{\ell_2=0}^{n_2}\sum_{m=\max\{\ell_1,\ell_2\}}^{n_1}A_{m,\ell_1,\ell_2}+\sum_{\ell_1=0}^{n_1}\sum_{\ell_2=0}^{n_2}\sum_{m=n_1+1}^{n_2}A_{m,\ell_1,\ell_2}+\sum_{\ell_1=0}^{n_1}\sum_{\ell_2=0}^{n_2}\sum_{m=n_2+1}^{\infty}A_{m,\ell_1,\ell_2},
\end{align*}
which coincides with Eq.~\eqref{sum3}.
\end{proof}
\begin{equation}\label{sum4}
\sum_{m=\max\{\ell_1,\ell_2\}}^{\infty}\sum_{k=\max\{0,\ell_1+\ell_2-m\}}^{\min\{\ell_1,\ell_2\}}A_{m,k}=\sum_{k=0}^{\min\{\ell_1,\ell_2\}}\sum_{m=\ell_1+\ell_2-k}^{\infty}A_{m,k}.
\end{equation}
\begin{proof}
Without the loss of generality, we assume that $\ell_1<\ell_2$.  Thus, Eq.~\eqref{sum4} is simplified into
\begin{equation}\label{sum4_1}
\sum_{m=\ell_2}^{\infty}\sum_{k=\max\{0,\ell_1+\ell_2-m\}}^{\ell_1}A_{m,k}=\sum_{k=0}^{\ell_1}\sum_{m=\ell_1+\ell_2-k}^{\infty}A_{m,k}.
\end{equation}
Then, the left-hand side of Eq.~\eqref{sum4_1} is split into
\begin{equation}\label{sum4_2}
\sum_{m=\ell_2}^{\infty}\sum_{k=\max\{0,\ell_1+\ell_2-m\}}^{\ell_1}A_{m,k}=\sum_{m=\ell_2}^{\ell_1+\ell_2}\sum_{k=\ell_1+\ell_2-m}^{\ell_1}A_{m,k}+\sum_{m=\ell_1+\ell_2+1}^{\infty}\sum_{k=0}^{\ell_1}A_{m,k}.
\end{equation}
In the right-hand side of Eq.~\eqref{sum4_2}, the sums of the second term are interchanged. The double summation of the first term is over the triangle $\ell_2\leq m\leq\ell_1+\ell_2$, $\ell_1+\ell_2-m\leq k\leq\ell_1$ which is rearranged into $0\leq k\leq\ell_1$, $\ell_1+\ell_2-k\leq m\leq \ell_1+\ell_2$. Thus:
\begin{equation*}
\sum_{m=\ell_2}^{\infty}\sum_{k=\max\{0,\ell_1+\ell_2-m\}}^{\ell_1}A_{m,k}=\sum_{k=0}^{\ell_1}\sum_{m=\ell_1+\ell_2-k}^{\ell_1+\ell_2}A_{m,k}+\sum_{k=0}^{\ell_1}\sum_{m=\ell_1+\ell_2+1}^{\infty}A_{m,k},
\end{equation*}
which coincides with Eq.~\eqref{sum4_1}.
\end{proof}

\section*{Acknowledgments}
The author is grateful to Prof. Jordan Stoyanov (Newcastle University/ Bulgarian Academy of Sciences) for his comments, which greatly improved the present work.  K.M. is supported by the Pacific Institute for the Mathematical Sciences (PIMS).

\bibliographystyle{plain}

\end{document}